

\input amstex

\def\filename{mrlt.sty}
\def\fileversion{2.1d}
\def\filedate{2-Sep-1999}
\expandafter\ifx\csname amsppt.sty\endcsname\endinput
  \expandafter\def\csname amsppt.sty\endcsname{2.1 (1-JUL-1991)}\fi
\xdef\fileversiontest{\fileversion\space(\filedate)}
\expandafter\ifx\csname\filename\endcsname\fileversiontest
  \message{[already loaded]}\endinput\fi
\expandafter\ifx\csname\filename\endcsname\relax 
  \else\errmessage{Discrepancy in `\filename' file versions:
     version \csname\filename\endcsname\space already loaded, trying
     now to load version \fileversiontest}\fi
\expandafter\xdef\csname\filename\endcsname{%
  \catcode`\noexpand\@=\the\catcode`\@
  \expandafter\gdef\csname\filename\endcsname{%
     \fileversion\space(\filedate)}}
\catcode`\@=11
\message{version \fileversion\space(\filedate):}
\expandafter\ifx\csname styname\endcsname\relax
  
\fi
\message{Loading utility definitions,}
\def\identity@#1{#1}
\def\nofrills@@#1{%
 \DN@{#1}%
 \ifx\next\nofrills \let\frills@\eat@
   \expandafter\expandafter\expandafter\next@\expandafter\eat@
  \else \let\frills@\identity@\expandafter\next@\fi}
\def\nofrillscheck#1{\def\nofrills@{\nofrills@@{#1}}%
  \futurelet\next\nofrills@}
\Invalid@\usualspace
\def\addto#1#2{\csname \expandafter\eat@\string#1@\endcsname
  \expandafter{\the\csname \expandafter\eat@\string#1@\endcsname#2}}
\newdimen\bigsize@
\def\big@#1#2{{\hbox{$\left#2\vcenter to#1\bigsize@{}%
  \right.\nulldelimiterspace\z@\m@th$}}}
\def\big{\big@\@ne}
\def\Big{\big@{1.5}}
\def\bigg{\big@\tw@}
\def\Bigg{\big@{2.5}}
\def\raggedcenter@{\leftskip\z@ plus.4\hsize \rightskip\leftskip
 \parfillskip\z@ \parindent\z@ \spaceskip.3333em \xspaceskip.5em
 \pretolerance9999\tolerance9999 \exhyphenpenalty\@M
 \hyphenpenalty\@M \let\\\linebreak}
\def\uppercasetext@#1{%
   {\spaceskip1.3\fontdimen2\the\font plus1.3\fontdimen3\the\font
    \def\ss{SS}\let\i=I\let\j=J\let\ae\AE\let\oe\OE
    \let\o\O\let\aa\AA\let\l\L
    \skipmath@#1$\skipmath@$}}
\def\skipmath@#1$#2${\uppercase{#1}%
  \ifx\skipmath@#2\else$#2$\expandafter\skipmath@\fi}
\def\add@missing#1{\expandafter\ifx\envir@end#1%
  \Err@{You seem to have a missing or misspelled
  \expandafter\string\envir@end ...}%
  \envir@end
\fi}
\newtoks\revert@
\def\envir@stack#1{\toks@\expandafter{\envir@end}%
  \edef\next@{\def\noexpand\envir@end{\the\toks@}%
    \revert@{\the\revert@}}%
  \revert@\expandafter{\next@}%
  \def\envir@end{#1}}
\begingroup
\catcode`\ =11
\gdef\revert@envir#1{\expandafter\ifx\envir@end#1%
\the\revert@%
\else\ifx\envir@end\enddocument \Err@{Extra \string#1}%
\else\expandafter\add@missing\envir@end\revert@envir#1%
\fi\fi}
\xdef\enddocument {\string\enddocument}%
\global\let\envir@end\enddocument 
\endgroup\relax
\def\first@#1#2\end{#1}
\def\true@{TT}
\def\false@{TF}
\def\empty@{}
\begingroup  \catcode`\-=3
\long\gdef\notempty#1{%
  \expandafter\ifx\first@#1-\end-\empty@ \false@\else \true@\fi}
\endgroup
\message{more fonts,}
\def\PSAMSFonts{TT}
\font@\tensmc=cmcsc10 \relax
\if\PSAMSFonts
  \font@\sevenex=cmex10 at 7pt
\else
  \font@\sevenex=cmex7 \relax
\fi
\font@\sevenit=cmti7 \relax
\font@\eightrm=cmr8 \relax 
\font@\sixrm=cmr6 \relax 
\font@\eighti=cmmi8 \relax     \skewchar\eighti='177 
\font@\sixi=cmmi6 \relax       \skewchar\sixi='177   
\font@\eightsy=cmsy8 \relax    \skewchar\eightsy='60 
\font@\sixsy=cmsy6 \relax      \skewchar\sixsy='60   
\if\PSAMSFonts
  \font@\eightex=cmex10 at 8pt
\else
  \font@\eightex=cmex8 \relax
\fi
\font@\eightbf=cmbx8 \relax 
\font@\sixbf=cmbx6 \relax   
\font@\eightit=cmti8 \relax 
\font@\eightsl=cmsl8 \relax 
\if\PSAMSFonts
  \font@\eightsmc=cmcsc10 at 8pt
\else
  \font@\eightsmc=cmcsc8 \relax
\fi
\font@\eighttt=cmtt8 \relax 
\loadeufm \loadmsam \loadmsbm
\message{symbol names}\UseAMSsymbols\message{,}
\newtoks\tenpoint@
\def\tenpoint{\normalbaselineskip12\p@
 \abovedisplayskip12\p@ plus3\p@ minus9\p@
 \belowdisplayskip\abovedisplayskip
 \abovedisplayshortskip\z@ plus3\p@
 \belowdisplayshortskip7\p@ plus3\p@ minus4\p@
 \textonlyfont@\rm\tenrm \textonlyfont@\it\tenit
 \textonlyfont@\sl\tensl \textonlyfont@\bf\tenbf
 \textonlyfont@\smc\tensmc \textonlyfont@\tt\tentt
 \ifsyntax@ \def\big##1{{\hbox{$\left##1\right.$}}}%
  \let\Big\big \let\bigg\big \let\Bigg\big
 \else
   \textfont\z@\tenrm  \scriptfont\z@\sevenrm
       \scriptscriptfont\z@\fiverm
   \textfont\@ne\teni  \scriptfont\@ne\seveni
       \scriptscriptfont\@ne\fivei
   \textfont\tw@\tensy \scriptfont\tw@\sevensy
       \scriptscriptfont\tw@\fivesy
   \textfont\thr@@\tenex \scriptfont\thr@@\sevenex
        \scriptscriptfont\thr@@\sevenex
   \textfont\itfam\tenit \scriptfont\itfam\sevenit
        \scriptscriptfont\itfam\sevenit
   \textfont\bffam\tenbf \scriptfont\bffam\sevenbf
        \scriptscriptfont\bffam\fivebf
   \setbox\strutbox\hbox{\vrule height8.5\p@ depth3.5\p@ width\z@}%
   \setbox\strutbox@\hbox{\lower.5\normallineskiplimit\vbox{%
        \kern-\normallineskiplimit\copy\strutbox}}%
   \setbox\z@\vbox{\hbox{$($}\kern\z@}\bigsize@1.2\ht\z@
  \fi
  \normalbaselines\rm\dotsspace@1.5mu\ex@.2326ex\jot3\ex@
  \the\tenpoint@}
\newtoks\eightpoint@
\def\eightpoint{\normalbaselineskip10\p@
 \abovedisplayskip10\p@ plus2.4\p@ minus7.2\p@
 \belowdisplayskip\abovedisplayskip
 \abovedisplayshortskip\z@ plus2.4\p@
 \belowdisplayshortskip5.6\p@ plus2.4\p@ minus3.2\p@
 \textonlyfont@\rm\eightrm \textonlyfont@\it\eightit
 \textonlyfont@\sl\eightsl \textonlyfont@\bf\eightbf
 \textonlyfont@\smc\eightsmc \textonlyfont@\tt\eighttt
 \ifsyntax@\def\big##1{{\hbox{$\left##1\right.$}}}%
  \let\Big\big \let\bigg\big \let\Bigg\big
 \else
  \textfont\z@\eightrm \scriptfont\z@\sixrm
       \scriptscriptfont\z@\fiverm
  \textfont\@ne\eighti \scriptfont\@ne\sixi
       \scriptscriptfont\@ne\fivei
  \textfont\tw@\eightsy \scriptfont\tw@\sixsy
       \scriptscriptfont\tw@\fivesy
  \textfont\thr@@\eightex \scriptfont\thr@@\sevenex
   \scriptscriptfont\thr@@\sevenex
  \textfont\itfam\eightit \scriptfont\itfam\sevenit
   \scriptscriptfont\itfam\sevenit
  \textfont\bffam\eightbf \scriptfont\bffam\sixbf
   \scriptscriptfont\bffam\fivebf
 \setbox\strutbox\hbox{\vrule height7\p@ depth3\p@ width\z@}%
 \setbox\strutbox@\hbox{\raise.5\normallineskiplimit\vbox{%
   \kern-\normallineskiplimit\copy\strutbox}}%
 \setbox\z@\vbox{\hbox{$($}\kern\z@}\bigsize@1.2\ht\z@
 \fi
 \normalbaselines\eightrm\dotsspace@1.5mu\ex@.2326ex\jot3\ex@
 \the\eightpoint@}
\def\linespacing#1{%
  \addto\tenpoint{\normalbaselineskip=#1\normalbaselineskip
    \normalbaselines
    \setbox\strutbox=\hbox{\vrule height.7\normalbaselineskip
      depth.3\normalbaselineskip}%
    \setbox\strutbox@\hbox{\raise.5\normallineskiplimit
      \vbox{\kern-\normallineskiplimit\copy\strutbox}}%
  }%
  \addto\eightpoint{\normalbaselineskip=#1\normalbaselineskip
    \normalbaselines
    \setbox\strutbox=\hbox{\vrule height.7\normalbaselineskip
      depth.3\normalbaselineskip}%
    \setbox\strutbox@\hbox{\raise.5\normallineskiplimit
      \vbox{\kern-\normallineskiplimit\copy\strutbox}}%
  }%
}
\if\PSAMSFonts
  \def\extrafont@#1#2#3{%
    \font#1=#2%
      \ifnum#3=9 10 at9pt%
      \else\ifnum#3=8 10 at8pt%
      \else\ifnum#3=6 7 at6pt%
              \else #3\fi\fi\fi\relax}
\else
  \def\extrafont@#1#2#3{\font#1=#2#3\relax}
\fi
\def\loadextrasizes@#1#2#3#4#5#6#7{%
 \ifx\undefined#1%
 \else \extrafont@{#4}{#2}{8}\extrafont@{#6}{#2}{6}%
   \ifsyntax@
   \else
     \addto\tenpoint{\textfont#1#3\scriptfont#1#5%
       \scriptscriptfont#1#7}%
    \addto\eightpoint{\textfont#1#4\scriptfont#1#6%
       \scriptscriptfont#1#7}%
   \fi
 \fi
}
\def\loadextrafonts@{%
  \loadextrasizes@\msafam{msam}%
    \tenmsa\eightmsa\sevenmsa\sixmsa\fivemsa
  \loadextrasizes@\msbfam{msbm}%
    \tenmsb\eightmsb\sevenmsb\sixmsb\fivemsb
  \loadextrasizes@\eufmfam{eufm}%
    \teneufm\eighteufm\seveneufm\sixeufm\fiveeufm
  \loadextrasizes@\eufbfam{eufb}%
    \teneufb\eighteufb\seveneufb\sixeufb\fiveeufb
  \loadextrasizes@\eusmfam{eusm}%
    \teneusm\eighteusm\seveneusm\sixeusm\fiveeusm
  \loadextrasizes@\eusbfam{eusb}%
    \teneusb\eighteusb\seveneusb\sixeusb\fiveeusb
  \loadextrasizes@\eurmfam{eurm}%
    \teneurm\eighteurm\seveneurm\sixeurm\fiveeurm
  \loadextrasizes@\eurbfam{eurb}%
    \teneurb\eighteurb\seveneurb\sixeurb\fiveeurb
  \loadextrasizes@\cmmibfam{cmmib}%
    \tencmmib\eightcmmib\sevencmmib\sixcmmib\fivecmmib
  \loadextrasizes@\cmbsyfam{cmbsy}%
    \tencmbsy\eightcmbsy\sevencmbsy\sixcmbsy\fivecmbsy
  \let\loadextrafonts@\empty@
}
\message{page dimension settings,}
\parindent1pc
\newdimen\normalparindent \normalparindent\parindent
\normallineskiplimit\p@
\newdimen\indenti \indenti=2pc
\def\pageheight#1{\vsize#1\relax}
\def\pagewidth#1{\hsize#1%
   \captionwidth@\hsize \advance\captionwidth@-2\indenti}
\pagewidth{30pc} \pageheight{47pc}
\let\magnification=\mag
\message{top matter,}
\def\topmatter{\loadextrafonts@ \let\topmatter\relax}
\def\chapterno@{\uppercase\expandafter{\romannumeral\chaptercount@}}
\newcount\chaptercount@
\def\chapter{\let\savedef@\chapter
  \def\chapter##1{\let\chapter\savedef@
  \leavevmode\hskip-\leftskip
   \rlap{\vbox to\z@{\vss\centerline{\eightpoint
   \frills@{CHAPTER\space\afterassignment\chapterno@
       \global\chaptercount@=}%
   ##1\unskip}\baselineskip2pc\null}}\hskip\leftskip}%
 \nofrillscheck\chapter}
\newbox\titlebox@

\def\title{\let\savedef@\title
 \def\title##1\endtitle{
   \let\title\savedef@
   \global\setbox\titlebox@\vtop{\tenpoint\bf
   \raggedcenter@
   \baselineskip1.3\baselineskip
   \frills@\uppercasetext@{##1}\endgraf}%
 \ifmonograph@ \edef\next{\the\leftheadtoks}%
    \ifx\next\empty@
    \leftheadtext{##1}\fi
 \fi
 \edef\next{\the\rightheadtoks}\ifx\next\empty@ \rightheadtext{##1}\fi
 }%
 \nofrillscheck\title}
\newbox\authorbox@
\def\author#1\endauthor{\global\setbox\authorbox@
 \vbox{\tenpoint\smc\raggedcenter@
 #1\endgraf}\relaxnext@ \edef\next{\the\leftheadtoks}%
 \ifx\next\empty@\leftheadtext{#1}\fi}
\newbox\affilbox@
\def\affil#1\endaffil{\global\setbox\affilbox@
 \vbox{\tenpoint\raggedcenter@#1\endgraf}}
\newcount\addresscount@
\addresscount@\z@
\def\address#1\endaddress{\global\advance\addresscount@\@ne
  \expandafter\gdef\csname address\number\addresscount@\endcsname
  {\nobreak\vskip12\p@ minus6\p@\indent\eightpoint\smc#1\par}}
\def\curraddr{\let\savedef@\curraddr
  \def\curraddr##1\endcurraddr{\let\curraddr\savedef@
  \toks@\expandafter\expandafter\expandafter{%
       \csname address\number\addresscount@\endcsname}%
  \toks@@{##1}%
  \expandafter\xdef\csname address\number\addresscount@\endcsname
  {\the\toks@\endgraf\noexpand\nobreak
    \indent{\noexpand\rm
    \frills@{{\noexpand\it Current address\noexpand\/}:\space}%
    \def\noexpand\usualspace{\space}\the\toks@@\unskip}}}%
  \nofrillscheck\curraddr}
\def\email{\let\savedef@\email
  \def\email##1\endemail{\let\email\savedef@
  \toks@{\def\usualspace{{\it\enspace}}\endgraf\indent\eightpoint}%
   \toks@@{\tt ##1\par}%
  \expandafter\xdef\csname email\number\addresscount@\endcsname
  {\the\toks@\frills@{{\noexpand\it E-mail address\noexpand\/}:%
     \noexpand\enspace}\the\toks@@}}%
  \nofrillscheck\email}
\def\thedate@{}
\def\date#1\enddate{\gdef\thedate@{\tenpoint#1\unskip}}
\def\thethanks@{}
\def\thanks#1\endthanks{%
  \ifx\thethanks@\empty@ \gdef\thethanks@{\eightpoint#1}%
  \else
    \expandafter\gdef\expandafter\thethanks@\expandafter{%
     \thethanks@\endgraf#1}%
  \fi}
\def\thekeywords@{}
\def\keywords{\let\savedef@\keywords
  \def\keywords##1\endkeywords{\let\keywords\savedef@
  \toks@{\def\usualspace{{\it\enspace}}\eightpoint}%
  \toks@@{##1\unskip.}%
  \edef\thekeywords@{\the\toks@\frills@{{\noexpand\it
    Key words and phrases.\noexpand\enspace}}\the\toks@@}}%
 \nofrillscheck\keywords}
\def\thesubjclass@{}
\def\subjclass{\let\savedef@\subjclass
 \def\subjclass##1\endsubjclass{\let\subjclass\savedef@
   \toks@{\def\usualspace{{\rm\enspace}}\eightpoint}%
   \toks@@{##1\unskip.}%
   \edef\thesubjclass@{\the\toks@
     \frills@{{\noexpand\rm2000 {\noexpand\it Mathematics Subject
       Classification}.\noexpand\enspace}}%
     \the\toks@@}}%
  \nofrillscheck\subjclass}
\newbox\abstractbox@
\def\abstract{\let\savedef@\abstract
 \def\abstract{\let\abstract\savedef@
  \setbox\abstractbox@\vbox\bgroup\noindent$$\vbox\bgroup
  \def\envir@end{\endabstract}\advance\hsize-2\indenti
  \def\usualspace{\enspace}\eightpoint \noindent
  \frills@{{\smc Abstract.\enspace}}}%
 \nofrillscheck\abstract}
\def\endabstract{\par\unskip\egroup$$\egroup}
\def\widestnumber{\begingroup \let\head\relax\let\subhead\relax
  \let\subsubhead\relax \expandafter\endgroup\setwidest@}
\def\setwidest@#1#2{%
   \ifx#1\head\setbox\tocheadbox@\hbox{#2.\enspace}%
   \else\ifx#1\subhead\setbox\tocsubheadbox@\hbox{#2.\enspace}%
   \else\ifx#1\subsubhead\setbox\tocsubheadbox@\hbox{#2.\enspace}%
   \else\ifx#1\key\refstyle A%
       \setboxz@h{\refsfont@\keyformat{#2}}%
       \refindentwd\wd\z@
   \else\ifx#1\no\refstyle C%
       \setboxz@h{\refsfont@\keyformat{#2}}%
       \refindentwd\wd\z@
   \else\ifx#1\page\setbox\z@\hbox{\quad\bf#2}%
       \pagenumwd\wd\z@
   \else\ifx#1\item
       \setboxz@h{(#2)}\rosteritemwd\wdz@
   \else\message{\string\widestnumber\space not defined for this
      option (\string#1)}%
\fi\fi\fi\fi\fi\fi\fi}
\newif\ifmonograph@
\def\Monograph{\monograph@true \let\headmark\rightheadtext
  \let\varindent@\indent \def\headfont@{\bf}\def\proclaimheadfont@{\smc}%
  \def\remarkheadfont@{\smc}}
\let\varindent@\noindent
\newbox\tocheadbox@    \newbox\tocsubheadbox@
\newbox\tocbox@
\newdimen\pagenumwd
\def\toc{\toc@{Contents}}
\def\newtocdefs{%
   \def \title##1\endtitle
       {\penaltyandskip@\z@\smallskipamount
        \hangindent\wd\tocheadbox@\noindent{\bf##1}}%
   \def \chapter##1{%
        Chapter \uppercase\expandafter{%
              \romannumeral##1.\unskip}\enspace}%
   \def \specialhead##1\endspecialhead
       {\par\hangindent\wd\tocheadbox@ \noindent##1\par}%
   \def \head##1 ##2\endhead
       {\par\hangindent\wd\tocheadbox@ \noindent
        \if\notempty{##1}\hbox to\wd\tocheadbox@{\hfil##1\enspace}\fi
        ##2\par}%
   \def \subhead##1 ##2\endsubhead
       {\par\vskip-\parskip {\normalbaselines
        \advance\leftskip\wd\tocheadbox@
        \hangindent\wd\tocsubheadbox@ \noindent
        \if\notempty{##1}%
              \hbox to\wd\tocsubheadbox@{##1\unskip\hfil}\fi
         ##2\par}}%
   \def \subsubhead##1 ##2\endsubsubhead
       {\par\vskip-\parskip {\normalbaselines
        \advance\leftskip\wd\tocheadbox@
        \hangindent\wd\tocsubheadbox@ \noindent
        \if\notempty{##1}%
              \hbox to\wd\tocsubheadbox@{##1\unskip\hfil}\fi
        ##2\par}}}
\def\toc@#1{\relaxnext@
 \DN@{\ifx\next\nofrills\DN@\nofrills{\nextii@}%
      \else\DN@{\nextii@{{#1}}}\fi
      \next@}%
 \DNii@##1{%
\ifmonograph@\bgroup\else\setbox\tocbox@\vbox\bgroup
   \centerline{\headfont@\ignorespaces##1\unskip}\nobreak
   \vskip\belowheadskip \fi
   \def\page####1%
       {\unskip\penalty\z@\null\hfil
        \rlap{\hbox to\pagenumwd{\quad\hfil####1}}%
              \hfilneg\penalty\@M}%
   \setbox\tocheadbox@\hbox{0.\enspace}%
   \setbox\tocsubheadbox@\hbox{0.0.\enspace}%
   \leftskip\indenti \rightskip\leftskip
   \setboxz@h{\bf\quad000}\pagenumwd\wd\z@
   \advance\rightskip\pagenumwd
   \newtocdefs
 }%
 \FN@\next@}
\def\endtoc{\par\egroup}
\let\pretitle\relax
\let\preauthor\relax
\let\preaffil\relax
\let\predate\relax
\let\preabstract\relax
\let\prepaper\relax
\def\dedicatory #1\enddedicatory{\def\preabstract{{\medskip
  \eightpoint\it \raggedcenter@#1\endgraf}}}
\def\thetranslator@{}
\def\translator{%
  \let\savedef@\translator
  \def\translator##1\endtranslator{\let\translator\savedef@
    \edef\thetranslator@{\noexpand\nobreak\noexpand\medskip
      \noexpand\line{\noexpand\eightpoint\hfil
      \frills@{Translated by \uppercase}{##1}\qquad\qquad}%
       \noexpand\nobreak}}%
  \nofrillscheck\translator}
\outer\def\endtopmatter{\add@missing\endabstract
 \edef\next{\the\leftheadtoks}\ifx\next\empty@
  \expandafter\leftheadtext\expandafter{\the\rightheadtoks}\fi
 \ifmonograph@\else
   \ifx\thesubjclass@\empty@\else \makefootnote@{}{\thesubjclass@}\fi
   \ifx\thekeywords@\empty@\else \makefootnote@{}{\thekeywords@}\fi
   \ifx\thethanks@\empty@\else \makefootnote@{}{\thethanks@}\fi
 \fi
  \pretitle
  \begingroup 
  \ifmonograph@ \topskip7pc \else \topskip4pc \fi
  \box\titlebox@
  \endgroup
  \preauthor
  \ifvoid\authorbox@\else \vskip2.5pcplus1pc\unvbox\authorbox@\fi
  \preaffil
  \ifvoid\affilbox@\else \vskip1pcplus.5pc\unvbox\affilbox@\fi
  \predate
  \ifx\thedate@\empty@\else
       \vskip1pcplus.5pc\line{\hfil\thedate@\hfil}\fi
  \preabstract
  \ifvoid\abstractbox@\else
       \vskip1.5pcplus.5pc\unvbox\abstractbox@ \fi
  \ifvoid\tocbox@\else\vskip1.5pcplus.5pc\unvbox\tocbox@\fi
  \prepaper
  \vskip2pcplus1pc\relax
}
\def\document{%
  \loadextrafonts@
  \let\fontlist@\relax\let\alloclist@\relax
  \tenpoint}
\message{section heads,}
\newskip\aboveheadskip       \aboveheadskip\bigskipamount
\newdimen\belowheadskip      \belowheadskip6\p@
\def\headfont@{\bf}
\def\penaltyandskip@#1#2{\par\skip@#2\relax
  \ifdim\lastskip<\skip@\relax\removelastskip
      \ifnum#1=\z@\else\penalty@#1\relax\fi\vskip\skip@
  \else\ifnum#1=\z@\else\penalty@#1\relax\fi\fi}
\def\nobreak{\penalty\@M
  \ifvmode\gdef\penalty@{\global\let\penalty@\penalty\count@@@}%
  \everypar{\global\let\penalty@\penalty\everypar{}}\fi}
\let\penalty@\penalty
\def\heading#1\endheading{\head#1\endhead}
\def\subheading{\DN@{\ifx\next\nofrills
    \expandafter\subheading@
  \else \expandafter\subheading@\expandafter\empty@
  \fi}%
  \FN@\next@
}
\def\subheading@#1#2{\subhead#1#2\endsubhead}
\def\specialheadfont@{\bf}
\outer\def\specialhead{%
  \add@missing\endroster \add@missing\enddefinition
  \add@missing\enddemo \add@missing\endexample
  \add@missing\endproclaim
  \penaltyandskip@{-200}\aboveheadskip
  \begingroup\interlinepenalty\@M\rightskip\z@ plus\hsize
  \let\\\linebreak
  \specialheadfont@\noindent}
\def\endspecialhead{\par\endgroup\nobreak\vskip\belowheadskip}
\outer\def\head#1\endhead{%
  \add@missing\endroster \add@missing\enddefinition
  \add@missing\enddemo \add@missing\endexample
  \add@missing\endproclaim
  \penaltyandskip@{-200}\aboveheadskip
  {\headfont@\raggedcenter@\interlinepenalty\@M
  #1\endgraf}\headmark{#1}%
  \nobreak
  \vskip\belowheadskip}
\let\headmark\eat@
\def\restoredef@#1{\relax\let#1\savedef@\let\savedef@\relax}
\newskip\subheadskip       \subheadskip\medskipamount
\def\subheadfont@{\bf}
\outer\def\subhead{%
  \add@missing\endroster \add@missing\enddefinition
  \add@missing\enddemo \add@missing\endexample
  \add@missing\endproclaim
  \let\savedef@\subhead \let\subhead\relax
  \def\subhead##1\endsubhead{\restoredef@\subhead
    \penaltyandskip@{-100}\subheadskip
    \varindent@{\def\usualspace{{\subheadfont@\enspace}}%
    \subheadfont@\ignorespaces##1\unskip\frills@{.\enspace}}%
    \ignorespaces}%
  \nofrillscheck\subhead}
\newskip\subsubheadskip       \subsubheadskip\medskipamount
\def\subsubheadfont@{\it}
\outer\def\subsubhead{%
  \add@missing\endroster \add@missing\enddefinition
  \add@missing\enddemo
  \add@missing\endexample \add@missing\endproclaim
  \let\savedef@\subsubhead \let\subsubhead\relax
  \def\subsubhead##1\endsubsubhead{\restoredef@\subsubhead
    \penaltyandskip@{-50}\subsubheadskip
      {\def\usualspace{\/{\it\enspace}}%
    \subsubheadfont@##1\unskip\frills@{.\enspace}}}%
  \nofrillscheck\subsubhead}
\message{theorems/proofs/definitions/remarks,}
\def\proclaimheadfont@{\bf}
\def\proclaimfont{\it}
\outer\def\proclaim{%
  \let\savedef@\proclaim \let\proclaim\relax
  \add@missing\endroster \add@missing\enddefinition
  \add@missing\endproclaim \envir@stack\endproclaim
 \def\proclaim##1{\restoredef@\proclaim
   \penaltyandskip@{-100}\medskipamount\varindent@
   \def\usualspace{{\proclaimheadfont@\enspace}}\proclaimheadfont@
   \ignorespaces##1\unskip\frills@{.\enspace}%
  \proclaimfont\ignorespaces}%
 \nofrillscheck\proclaim}
\def\endproclaim{\revert@envir\endproclaim \par\rm
  \penaltyandskip@{55}\medskipamount}
\def\remarkheadfont@{\it}
\def\remark{\let\savedef@\remark \let\remark\relax
  \add@missing\endroster \add@missing\endproclaim
  \envir@stack\endremark
  \def\remark##1{\restoredef@\remark
    \penaltyandskip@\z@\medskipamount
  {\def\usualspace{{\remarkheadfont@\enspace}}%
  \varindent@\remarkheadfont@\ignorespaces##1\unskip%
  \frills@{.\enspace}}\rm
  \ignorespaces}\nofrillscheck\remark}
\def\endremark{\par\revert@envir\endremark}
\ifx\undefined\square
  \def\square{\vrule width.6em height.5em depth.1em\relax}\fi
\def\qed{\ifhmode\unskip\nobreak\fi\quad
  \ifmmode\square\else$\m@th\square$\fi}
\def\demo{\DN@{\ifx\next\nofrills
    \DN@####1####2{\remark####1{####2}\envir@stack\enddemo
      \ignorespaces}%
  \else
    \DN@####1{\remark{####1}\envir@stack\enddemo\ignorespaces}%
  \fi
  \next@}%
\FN@\next@}

\def\enddemo{\par\revert@envir\enddemo \endremark\medskip}
\def\definition{\let\savedef@\definition \let\definition\relax
  \add@missing\endproclaim \add@missing\endroster
  \add@missing\enddefinition \envir@stack\enddefinition
   \def\definition##1{\restoredef@\definition
     \penaltyandskip@{-100}\medskipamount
        {\def\usualspace{{\proclaimheadfont@\enspace}}%
        \varindent@\proclaimheadfont@\ignorespaces##1\unskip
        \frills@{.\proclaimheadfont@\enspace}}%
        \rm \ignorespaces}%
  \nofrillscheck\definition}
\def\enddefinition{\revert@envir\enddefinition
  \par\medskip}
\def\example{\DN@{\ifx\next\nofrills
    \DN@####1####2{\definition####1{####2}\envir@stack\endexample
      \ignorespaces}%
  \else
    \DN@####1{\definition{####1}\envir@stack\endexample\ignorespaces}%
  \fi
  \next@}%
\FN@\next@}
\def\endexample{\revert@envir\endexample \enddefinition }
\message{rosters,}
\newdimen\rosteritemwd
\rosteritemwd16pt 
\newcount\rostercount@
\newif\iffirstitem@
\let\plainitem@\item
\newtoks\everypartoks@
\def\par@{\everypartoks@\expandafter{\the\everypar}\everypar{}}
\def\leftskip@{}
\def\roster{%
  \envir@stack\endroster
 \edef\leftskip@{\leftskip\the\leftskip}%
 \relaxnext@
 \rostercount@\z@
 \def\item{\FN@\rosteritem@}
 \DN@{\ifx\next\runinitem\let\next@\nextii@\else
  \let\next@\nextiii@\fi\next@}%
 \DNii@\runinitem
  {\unskip
   \DN@{\ifx\next[\let\next@\nextii@\else
    \ifx\next"\let\next@\nextiii@\else\let\next@\nextiv@\fi\fi\next@}%
   \DNii@[####1]{\rostercount@####1\relax
    \enspace\therosteritem{\number\rostercount@}~\ignorespaces}%
   \def\nextiii@"####1"{\enspace{\rm####1}~\ignorespaces}%
   \def\nextiv@{\enspace\therosteritem1\rostercount@\@ne~}%
   \par@\firstitem@false
   \FN@\next@}
 \def\nextiii@{\par\par@
  \penalty\@m\smallskip\vskip-\parskip
  \firstitem@true}
 \FN@\next@}
\def\rosteritem@{\iffirstitem@\firstitem@false
  \else\par\vskip-\parskip\fi
 \leftskip\rosteritemwd \advance\leftskip\normalparindent
 \advance\leftskip.5em \noindent
 \DNii@[##1]{\rostercount@##1\relax\itembox@}%
 \def\nextiii@"##1"{\def\therosteritem@{\rm##1}\itembox@}%
 \def\nextiv@{\advance\rostercount@\@ne\itembox@}%
 \def\therosteritem@{\therosteritem{\number\rostercount@}}%
 \ifx\next[\let\next@\nextii@\else\ifx\next"\let\next@\nextiii@\else
  \let\next@\nextiv@\fi\fi\next@}
\def\itembox@{\llap{\hbox to\rosteritemwd{\hss
  \kern\z@ 
  \therosteritem@}\enspace}\ignorespaces}
\def\therosteritem#1{\rom{(\ignorespaces#1\unskip)}}
\newif\ifnextRunin@
\def\endroster{\relaxnext@
 \revert@envir\endroster 
 \par\leftskip@
 \global\rosteritemwd16\p@ 
 \penalty-50 \vskip-\parskip\smallskip
 \DN@{\ifx\next\Runinitem\let\next@\relax
  \else\nextRunin@false\let\item\plainitem@
   \ifx\next\par
    \DN@\par{\everypar\expandafter{\the\everypartoks@}}%
   \else
    \DN@{\noindent\everypar\expandafter{\the\everypartoks@}}%
  \fi\fi\next@}%
 \FN@\next@}
\newcount\rosterhangafter@
\def\Runinitem#1\roster\runinitem{\relaxnext@
  \envir@stack\endroster
 \rostercount@\z@
 \def\item{\FN@\rosteritem@}%
 \def\runinitem@{#1}%
 \DN@{\ifx\next[\let\next\nextii@\else\ifx\next"\let\next\nextiii@
  \else\let\next\nextiv@\fi\fi\next}%
 \DNii@[##1]{\rostercount@##1\relax
  \def\item@{\therosteritem{\number\rostercount@}}\nextv@}%
 \def\nextiii@"##1"{\def\item@{{\rm##1}}\nextv@}%
 \def\nextiv@{\advance\rostercount@\@ne
  \def\item@{\therosteritem{\number\rostercount@}}\nextv@}%
 \def\nextv@{\setbox\z@\vbox
  {\ifnextRunin@\noindent\fi
  \runinitem@\unskip\enspace\item@~\par
  \global\rosterhangafter@\prevgraf}%
  \firstitem@false
  \ifnextRunin@\else\par\fi
  \hangafter\rosterhangafter@\hangindent3\normalparindent
  \ifnextRunin@\noindent\fi
  \runinitem@\unskip\enspace
  \item@~\ifnextRunin@\else\par@\fi
  \nextRunin@true\ignorespaces}
 \FN@\next@}
\message{footnotes,}
\def\footmarkform@#1{$\m@th^{#1}$}
\let\thefootnotemark\footmarkform@
\def\makefootnote@#1#2{\insert\footins
 {\interlinepenalty\interfootnotelinepenalty
 \eightpoint\splittopskip\ht\strutbox\splitmaxdepth\dp\strutbox
 \floatingpenalty\@MM\leftskip\z@skip\rightskip\z@skip
 \spaceskip\z@skip\xspaceskip\z@skip
 \leavevmode{#1}\footstrut\ignorespaces#2\unskip\lower\dp\strutbox
 \vbox to\dp\strutbox{}}}
\newcount\footmarkcount@
\footmarkcount@\z@
\def\footnotemark{\let\@sf\empty@\relaxnext@
 \ifhmode\edef\@sf{\spacefactor\the\spacefactor}\/\fi
 \DN@{\ifx[\next\let\next@\nextii@\else
  \ifx"\next\let\next@\nextiii@\else
  \let\next@\nextiv@\fi\fi\next@}%
 \DNii@[##1]{\footmarkform@{##1}\@sf}%
 \def\nextiii@"##1"{{##1}\@sf}%
 \def\nextiv@{\iffirstchoice@\global\advance\footmarkcount@\@ne\fi
  \footmarkform@{\number\footmarkcount@}\@sf}%
 \FN@\next@}
\def\footnotetext{\relaxnext@
 \DN@{\ifx[\next\let\next@\nextii@\else
  \ifx"\next\let\next@\nextiii@\else
  \let\next@\nextiv@\fi\fi\next@}%
 \DNii@[##1]##2{\makefootnote@{\footmarkform@{##1}}{##2}}%
 \def\nextiii@"##1"##2{\makefootnote@{##1}{##2}}%
 \def\nextiv@##1{\makefootnote@{\footmarkform@%
  {\number\footmarkcount@}}{##1}}%
 \FN@\next@}
\def\footnote{\let\@sf\empty@\relaxnext@
 \ifhmode\edef\@sf{\spacefactor\the\spacefactor}\/\fi
 \DN@{\ifx[\next\let\next@\nextii@\else
  \ifx"\next\let\next@\nextiii@\else
  \let\next@\nextiv@\fi\fi\next@}%
 \DNii@[##1]##2{\footnotemark[##1]\footnotetext[##1]{##2}}%
 \def\nextiii@"##1"##2{\footnotemark"##1"\footnotetext"##1"{##2}}%
 \def\nextiv@##1{\footnotemark\footnotetext{##1}}%
 \FN@\next@}
\def\adjustfootnotemark#1{\advance\footmarkcount@#1\relax}
\def\footnoterule{\kern-3\p@
  \hrule width5pc\kern 2.6\p@}
\message{figures and captions,}
\def\captionfont@{\smc}
\def\topcaption#1#2\endcaption{%
  {\dimen@\hsize \advance\dimen@-\captionwidth@
   \rm\raggedcenter@ \advance\leftskip.5\dimen@ \rightskip\leftskip
  {\captionfont@#1}%
  \if\notempty{#2}.\enspace\ignorespaces#2\fi
  \endgraf}\nobreak\bigskip}
\def\botcaption#1#2\endcaption{%
  \nobreak\bigskip
  \setboxz@h{\captionfont@#1\if\notempty{#2}.\enspace\rm#2\fi}%
  {\dimen@\hsize \advance\dimen@-\captionwidth@
   \leftskip.5\dimen@ \rightskip\leftskip
   \noindent \ifdim\wdz@>\captionwidth@
   \else\hfil\fi
  {\captionfont@#1}%
  \if\notempty{#2}.\enspace\rm#2\fi\endgraf}}
\def\@ins{\par\begingroup\def\vspace##1{\vskip##1\relax}%
  \def\captionwidth##1{\captionwidth@##1\relax}%
  \setbox\z@\vbox\bgroup} 
\message{miscellaneous,}
\def\block{\RIfMIfI@\nondmatherr@\block\fi
       \else\ifvmode\noindent$$\predisplaysize\hsize
         \else$$\fi
  \def\endblock{\par\egroup$$}\fi
  \vbox\bgroup\advance\hsize-2\indenti\noindent}
\def\endblock{\par\egroup}
\def\cite#1{\rom{[{\citefont@\m@th#1}]}}
\def\citefont@{\rm}
\def\rom#1{\leavevmode
  \edef\prevskip@{\ifdim\lastskip=\z@ \else\hskip\the\lastskip\relax\fi}%
  \unskip
  \edef\prevpenalty@{\ifnum\lastpenalty=\z@ \else
    \penalty\the\lastpenalty\relax\fi}%
  \unpenalty \/\prevpenalty@ \prevskip@ {\rm #1}}
\message{references,}
\def\refsfont@{\eightpoint}
\newdimen\refindentwd
\setboxz@h{\refsfont@ 00.\enspace}
\refindentwd\wdz@
\outer\def\Refs{\add@missing\endroster \add@missing\endproclaim
 \let\savedef@\Refs \let\Refs\relax 
 \def\Refs##1{\restoredef@\Refs
   \if\notempty{##1}\penaltyandskip@{-200}\aboveheadskip
     \begingroup \raggedcenter@\headfont@
       \ignorespaces##1\endgraf\endgroup
     \penaltyandskip@\@M\belowheadskip
   \fi
   \begingroup\def\envir@end{\endRefs}\refsfont@\sfcode`\.\@m
   }%
 \nofrillscheck{\csname Refs\expandafter\endcsname
  \frills@{{References}}}}
\def\endRefs{\par 
  \endgroup}
\newif\ifbook@ \newif\ifprocpaper@
\def\nofrills{%
  \expandafter\ifx\envir@end\endref
    \let\do\relax
    \xdef\nofrills@list{\nofrills@list\do\curbox}%
  \else\errmessage{\Invalid@@ \string\nofrills}%
  \fi}%
\def\defaultreftexts{\gdef\edtext{ed.}\gdef\pagestext{pp.}%
  \gdef\voltext{vol.}\gdef\issuetext{no.}}
\defaultreftexts
\def\ref{\par
  \begingroup \def\envir@end{\endref}%
  \noindent\hangindent\refindentwd
  \def\par{\add@missing\endref}%
  \global\let\nofrills@list\empty@
  \refbreaks
  \procpaper@false \book@false \moreref@false
  \def\curbox{\z@}\setbox\z@\vbox\bgroup
}
\let\keyhook@\empty@
\def\endref{%
  \setbox\tw@\box\thr@@
  \makerefbox?\thr@@{\endgraf\egroup}%
  \endref@
  \endgraf
  \endgroup
  \keyhook@
  \global\let\keyhook@\empty@ 
}
\def\key{\gdef\key{\makerefbox\key\keybox@\empty@}\key} \newbox\keybox@
\def\no{\gdef\no{\makerefbox\no\keybox@\empty@}%
  \gdef\keyhook@{\refstyle C}\no}
\def\by{\makerefbox\by\bybox@\empty@} \newbox\bybox@
\def\bysame{\by\hbox to3em{\hrulefill}\thinspace\kern\z@}
\def\paper{\makerefbox\paper\paperbox@\it} \newbox\paperbox@
\def\paperinfo{\makerefbox\paperinfo\paperinfobox@\empty@}%
  \newbox\paperinfobox@
\def\jour{\makerefbox\jour\jourbox@
  {\aftergroup\book@false \aftergroup\procpaper@false}} \newbox\jourbox@
\def\issue{\makerefbox\issue\issuebox@\empty@} \newbox\issuebox@
\def\yr{\makerefbox\yr\yrbox@\empty@} \newbox\yrbox@
\def\pages{\makerefbox\pages\pagesbox@\empty@} \newbox\pagesbox@
\def\page{\gdef\pagestext{p.}\makerefbox\page\pagesbox@\empty@}
\def\ed{\makerefbox\ed\edbox@\empty@} \newbox\edbox@
\def\eds{\gdef\edtext{eds.}\makerefbox\eds\edbox@\empty@}
\def\book{\makerefbox\book\bookbox@
  {\it\aftergroup\book@true \aftergroup\procpaper@false}}
  \newbox\bookbox@
\def\bookinfo{\makerefbox\bookinfo\bookinfobox@\empty@}%
  \newbox\bookinfobox@
\def\publ{\makerefbox\publ\publbox@\empty@} \newbox\publbox@
\def\publaddr{\makerefbox\publaddr\publaddrbox@\empty@}%
  \newbox\publaddrbox@
\def\inbook{\makerefbox\inbook\bookbox@
  {\aftergroup\procpaper@true \aftergroup\book@false}}
\def\procinfo{\makerefbox\procinfo\procinfobox@\empty@}%
  \newbox\procinfobox@
\def\finalinfo{\makerefbox\finalinfo\finalinfobox@\empty@}%
  \newbox\finalinfobox@
\def\miscnote{\makerefbox\miscnote\miscnotebox@\empty@}%
  \newbox\miscnotebox@

\def\lang{\makerefbox\lang\langbox@\empty@} \newbox\langbox@
\newbox\morerefbox@
\def\vol{\makerefbox\vol\volbox@{\ifbook@ \else
  \ifprocpaper@\else\bf\fi\fi}}
\newbox\volbox@
\newbox\holdoverbox
\def\makerefbox#1#2#3{\endgraf
  \setbox\z@\lastbox
  \global\setbox\@ne\hbox{\unhbox\holdoverbox
    \ifvoid\z@\else\unhbox\z@\unskip\unskip\unpenalty\fi}%
  \egroup
  \setbox\curbox\box\ifdim\wd\@ne>\z@ \@ne \else\voidb@x\fi
  \ifvoid#2\else\Err@{Redundant \string#1; duplicate use, or
     mutually exclusive information already given}\fi
  \def\curbox{#2}\setbox\curbox\vbox\bgroup \hsize\maxdimen \noindent
  #3}
\def\refbreaks{%
  \def\refconcat##1{\setbox\z@\lastbox \setbox\holdoverbox\hbox{%
       \unhbox\holdoverbox \unhbox\z@\unskip\unskip\unpenalty##1}}%
  \def\holdover##1{%
    \RIfM@
      \penalty-\@M\null
      \hfil$\clubpenalty\z@\widowpenalty\z@\interlinepenalty\z@
      \offinterlineskip\endgraf
      \setbox\z@\lastbox\unskip \unpenalty
      \refconcat{##1}%
      \noindent
      $\hfil\penalty-\@M
    \else
      \endgraf\refconcat{##1}\noindent
    \fi}%
  \def\break{\holdover{\penalty-\@M}}%
  \let\vadjust@\vadjust
  \def\vadjust##1{\holdover{\vadjust@{##1}}}%
  \def\newpage{\vadjust{\vfill\break}}%
}
\def\refstyle#1{\uppercase{%
  \if#1A\relax \def\keyformat##1{[##1]\enspace\hfil}%
  \else\if#1B\relax
    \def\keyformat##1{\aftergroup\kern
              \aftergroup-\aftergroup\refindentwd}%
    \refindentwd\parindent
 \else\if#1C\relax
   \def\keyformat##1{\hfil##1.\enspace}%
 \fi\fi\fi}
}
\refstyle{A}
\def\finalpunct{\ifnum\lastkern=\m@ne\unkern\else.\fi
       \refquotes@\refbreak@}%
\def\continuepunct#1#2#3#4{}%
\def\endref@{%
  \keyhook@
  \def\nofrillscheck##1{%
    \def\do####1{\ifx##1####1\let\frills@\eat@\fi}%
    \let\frills@\identity@ \nofrills@list}%
  \ifvoid\bybox@
    \ifvoid\edbox@
    \else\setbox\bybox@\hbox{\unhbox\edbox@\breakcheck
      \nofrillscheck\edbox@\frills@{\space(\edtext)}\refbreak@}\fi
  \fi
  \ifvoid\keybox@\else\hbox to\refindentwd{%
       \keyformat{\unhbox\keybox@}}\fi
  \ifmoreref@
    \commaunbox@\morerefbox@
  \else
    \kern-\tw@ sp\kern\m@ne sp
  \fi
  \ppunbox@\empty@\empty@\bybox@\empty@
  \ifbook@ 
    \commaunbox@\bookbox@ \commaunbox@\bookinfobox@
    \ppunbox@\empty@{ (}\procinfobox@)%
    \ppunbox@,{ vol.~}\volbox@\empty@
    \ppunbox@\empty@{ (}\edbox@{, \edtext)}%
    \commaunbox@\publbox@ \commaunbox@\publaddrbox@
    \commaunbox@\yrbox@
    \ppunbox@,{ \pagestext~}\pagesbox@\empty@
  \else
    \commaunbox@\paperbox@ \commaunbox@\paperinfobox@
    \ifprocpaper@ 
      \commaunbox@\bookbox@
      \ppunbox@\empty@{ (}\procinfobox@)%
      \ppunbox@\empty@{ (}\edbox@{, \edtext)}%
      \commaunbox@\bookinfobox@
      \ppunbox@,{ \voltext~}\volbox@\empty@
      \commaunbox@\publbox@ \commaunbox@\publaddrbox@
      \commaunbox@\yrbox@
      \ppunbox@,{ \pagestext~}\pagesbox@\empty@
    \else 
      \commaunbox@\jourbox@
      \ppunbox@\empty@{ }\volbox@\empty@
      \ppunbox@\empty@{ (}\yrbox@)%
      \ppunbox@,{ \issuetext~}\issuebox@\empty@
      \commaunbox@\publbox@ \commaunbox@\publaddrbox@
      \commaunbox@\pagesbox@
    \fi
  \fi
  \commaunbox@\finalinfobox@
  \ppunbox@\empty@{ (}\miscnotebox@)%
  \finalpunct\ppunbox@\empty@{ (}\langbox@)%
  \defaultreftexts
}
\def\punct@#1{#1}
\def\ppunbox@#1#2#3#4{\ifvoid#3\else
  \let\prespace@\relax
  \ifnum\lastkern=\m@ne \unkern\let\punct@\eat@
    \ifnum\lastkern=-\tw@ \unkern\let\prespace@\ignorespaces \fi
  \fi
  \nofrillscheck#3%
  \punct@{#1}\refquotes@\refbreak@
  \let\punct@\identity@
  \prespace@
  \frills@{#2\eat@}\space
  \unhbox#3\breakcheck
  \frills@{#4\eat@}{\kern\m@ne sp}\fi}
\def\commaunbox@#1{\ppunbox@,\space{#1}\empty@}
\def\breakcheck{\edef\refbreak@{\ifnum\lastpenalty=\z@\else
  \penalty\the\lastpenalty\relax\fi}\unpenalty}
\def\endquotes{\def\refquotes@{''\let\refquotes@\empty@}}
\let\refquotes@\empty@
\let\refbreak@\empty@
\newif\ifmoreref@
\def\moreref{%
  \setbox\tw@\box\thr@@
  \makerefbox?\thr@@{\endgraf\egroup}%
  \let\savedef@\finalpunct  \let\finalpunct\empty@
  \endref@
  \def\punct@##1##2{##2;}%
  \global\let\nofrills@list\empty@ 
  \let\finalpunct\savedef@
  \moreref@true
  \def\curbox{\morerefbox@}%
  \setbox\morerefbox@\vbox\bgroup \hsize\maxdimen \noindent
}

\message{end of document,}
\outer\def\enddocument{\par
  \add@missing\endRefs
  \add@missing\endroster \add@missing\endproclaim
  \add@missing\enddefinition
  \add@missing\enddemo \add@missing\endremark \add@missing\endexample
 \ifmonograph@ 
 \else
 \nobreak
 \thetranslator@
 \count@\z@ \loop\ifnum\count@<\addresscount@\advance\count@\@ne
 \csname address\number\count@\endcsname
 \csname email\number\count@\endcsname
 \repeat
\fi
 \vfill\supereject\end}
\message{output routine,}
\def\folio{{\foliofont@\ifnum\pageno<\z@ \romannumeral-\pageno
 \else\number\pageno \fi}}
\def\foliofont@{\eightrm}
\def\headlinefont@{\eightpoint}
\def\leftheadline{\rlap{\folio}\hfill \iftrue\topmark\fi \hfill}
\def\rightheadline{\hfill \expandafter
  \hfill \llap{\folio}}
\newtoks\leftheadtoks
\newtoks\rightheadtoks
\def\leftheadtext{\let\savedef@\leftheadtext
  \def\leftheadtext##1{\let\leftheadtext\savedef@
    \leftheadtoks\expandafter{\frills@\uppercasetext@{##1}}%
    \mark{\the\leftheadtoks\noexpand\else\the\rightheadtoks}
    \ifsyntax@\setboxz@h{\def\\{\unskip\space\ignorespaces}%
        \headlinefont@##1}\fi}%
  \nofrillscheck\leftheadtext}
\def\rightheadtext{\let\savedef@\rightheadtext
  \def\rightheadtext##1{\let\rightheadtext\savedef@
    \rightheadtoks\expandafter{\frills@\uppercasetext@{##1}}%
    \mark{\the\leftheadtoks\noexpand\else\the\rightheadtoks}%
    \ifsyntax@\setboxz@h{\def\\{\unskip\space\ignorespaces}%
        \headlinefont@##1}\fi}%
  \nofrillscheck\rightheadtext}
\headline={\def\\{\unskip\space\ignorespaces}\headlinefont@
  \def\chapter{%
    \def\chapter##1{%
      \frills@{\afterassignment\chapterno@ \chaptercount@=}##1.\space}%
    \nofrillscheck\chapter}%
  \ifodd\pageno \rightheadline \else \leftheadline\fi}
\def\NoRunningHeads{\global\runheads@false\global\let\headmark\eat@}

\def\logo@{\baselineskip2pc \hbox to\hsize{\hfil\eightpoint Typeset by
 \AmSTeX}}
\def\nologo{\def\logo@{}}
\newif\iffirstpage@     \firstpage@true
\newif\ifrunheads@      \runheads@true
\output={\output@}
\def\output@{\shipout\vbox{%
 \iffirstpage@ \global\firstpage@false
  \pagebody \logo@ \makefootline%
 \else \ifrunheads@ \makeheadline \pagebody
       \else \pagebody \makefootline \fi
 \fi}%
 \advancepageno \ifnum\outputpenalty>-\@MM\else\dosupereject\fi}
\message{hyphenation exceptions (U.S. English)}
\hyphenation{acad-e-my acad-e-mies af-ter-thought anom-aly anom-alies
an-ti-deriv-a-tive an-tin-o-my an-tin-o-mies apoth-e-o-ses
apoth-e-o-sis ap-pen-dix ar-che-typ-al as-sign-a-ble as-sist-ant-ship
as-ymp-tot-ic asyn-chro-nous at-trib-uted at-trib-ut-able bank-rupt
bank-rupt-cy bi-dif-fer-en-tial blue-print busier busiest
cat-a-stroph-ic cat-a-stroph-i-cally con-gress cross-hatched data-base
de-fin-i-tive de-riv-a-tive dis-trib-ute dri-ver dri-vers eco-nom-ics
econ-o-mist elit-ist equi-vari-ant ex-quis-ite ex-tra-or-di-nary
flow-chart for-mi-da-ble forth-right friv-o-lous ge-o-des-ic
ge-o-det-ic geo-met-ric griev-ance griev-ous griev-ous-ly
hexa-dec-i-mal ho-lo-no-my ho-mo-thetic ideals idio-syn-crasy
in-fin-ite-ly in-fin-i-tes-i-mal ir-rev-o-ca-ble key-stroke
lam-en-ta-ble light-weight mal-a-prop-ism man-u-script mar-gin-al
meta-bol-ic me-tab-o-lism meta-lan-guage me-trop-o-lis
met-ro-pol-i-tan mi-nut-est mol-e-cule mono-chrome mono-pole
mo-nop-oly mono-spline mo-not-o-nous mul-ti-fac-eted mul-ti-plic-able
non-euclid-ean non-iso-mor-phic non-smooth par-a-digm par-a-bol-ic
pa-rab-o-loid pa-ram-e-trize para-mount pen-ta-gon phe-nom-e-non
post-script pre-am-ble pro-ce-dur-al pro-hib-i-tive pro-hib-i-tive-ly
pseu-do-dif-fer-en-tial pseu-do-fi-nite pseu-do-nym qua-drat-ic
quad-ra-ture qua-si-smooth qua-si-sta-tion-ary qua-si-tri-an-gu-lar
quin-tes-sence quin-tes-sen-tial re-arrange-ment rec-tan-gle
ret-ri-bu-tion retro-fit retro-fit-ted right-eous right-eous-ness
ro-bot ro-bot-ics sched-ul-ing se-mes-ter semi-def-i-nite
semi-ho-mo-thet-ic set-up se-vere-ly side-step sov-er-eign spe-cious
spher-oid spher-oid-al star-tling star-tling-ly sta-tis-tics
sto-chas-tic straight-est strange-ness strat-a-gem strong-hold
sum-ma-ble symp-to-matic syn-chro-nous topo-graph-i-cal tra-vers-a-ble
tra-ver-sal tra-ver-sals treach-ery turn-around un-at-tached
un-err-ing-ly white-space wide-spread wing-spread wretch-ed
wretch-ed-ly Brown-ian Eng-lish Euler-ian Feb-ru-ary Gauss-ian
Grothen-dieck Hamil-ton-ian Her-mit-ian Jan-u-ary Japan-ese Kor-te-weg
Le-gendre Lip-schitz Lip-schitz-ian Mar-kov-ian Noe-ther-ian
No-vem-ber Rie-mann-ian Schwarz-schild Sep-tem-ber}
\tenpoint
\W@{}
\csname mrlt.sty\endcsname

\countdef\pageno=0

\baselineskip=13pt
\hsize=122mm
\vsize=187mm
\magnification=1100

\nologo

\def\section#1{\heading{#1}\endheading}
\def\subsection#1{\subheading{#1}}




\pageno=1



\def\rightheadline{\eightpoint\kern3.5em
IN DIMENSION 2, COMPLETE IDEALS ARE MULTIPLIER IDEALS\kern3.1em\the\pageno}


\nologo

\def\proof{\demo{Proof}}
\def\endproof{\qed\enddemo}
\def\nextpart#1{\medskip{\bf (#1).}}

\def\>{\mkern1mu}
\def\<{\mkern-1mu}
\def\set{\!:=}
\def\nmb{\nomathbreak}
\def\smcirc{{\hbox to.7em{$\hss \eightpoint\circ \hss$}}}
\def\blp{\text{\bf(}}
\def\brp{\text{\bf)}}
\def\iso
   {{\mkern8mu\longrightarrow \mkern-25.5mu{}^\sim\mkern17mu}}

\def\Spec{\text{\rm Spec}}
\def\Div{\operatorname{Div}}
\def\Dive{\operatorname{Div}_{\text{\rm e}}}
\def\J{{\Cal J}}
\def\L{{\Cal L}}
\def\O{{\Cal O}}
\def\OD{{\O_{\mkern-2.5mu D}}}
\def\OE{{\O_{\mkern-2.5mu E}}}
\def\OX{{\O_{\mkern-2.5mu X}}}
\def\aff{{_{\<\text{aff}}}}
\def\exc{{_{\<\text{exc}}}}

\loadbold


\topmatter

\subjclass  13B22, 13H05 \endsubjclass

\title {Integrally closed ideals\\ in two-dimensional 
regular local rings\\ are multiplier ideals}
\endtitle

\author Joseph Lipman and Kei-ichi Watanabe\endauthor


\address{Dept\. of Mathematics, Purdue University,
 W. Lafayette IN 47907, USA} \endaddress

\email{lipman\@math.purdue.edu}\endemail


\address{Dept\. of Mathematics, Nihon University, Sakura Josui 3-25-40,
Setagaya, Tokyo 156-8550, Japan  } \endaddress

\email{watanabe\@math.chs.nihon-u.ac.jp}   \endemail


\thanks First author partially supported by the National Security
Agency. Second author partially supported by Grants-in-Aid
in Scientific Researches, 13440015,  13874006; and his stay at MSRI was
supported by the Bunri Fund, Nihon University.
Both authors are grateful to MSRI for providing the environment
without which this work would not have begun.
Research at MSRI is supported in part by NSF grant DMS-9810361.
\vadjust{\kern1pt} 
\endthanks

\abstract{
Multiplier ideals in commutative rings are certain integrally closed ideals 
with properties that lend themselves to highly interesting applications. 
How special are they among integrally closed ideals in general? 
We show that in a two-dimensional regular local ring with 
algebraically closed residue field  there is in fact no
difference between ``multiplier" and ``integrally closed" (or \hbox{``complete."})
But among  multiplier ideals arising from an {\it integer\/} multiplying constant (also known
as  {\it adjoint\/} ideals), and primary for the maximal ideal, the only simple complete ideals
are those of order one.}\looseness=-1
\endabstract

\endtopmatter


\document

\subheading{Introduction} There has arisen in recent years a substantial
body of work on {\it multiplier ideals} in commutative rings (see \cite{La}). 
Multiplier ideals are integrally closed ideals with properties that lend
themselves to highly interesting applications. One is tempted then to
ask just how special multiplier ideals are
among integrally closed ideals in general. 

In this note we show that in a two-dimensional regular local ring $R$
with maximal ideal $\frak m$ such that the residue field $R/\frak m$
is algebraically closed%
\footnote{It most likely suffices that $R/\frak m$ be infinite, but we 
want to avoid additional technicalities.} 
there is actually no difference between
multiplier ideals and integrally closed ideals. In fact it turns out more convenient to do this
for {\it fractionary\/} $R$-ideals, i.e., nonzero finitely-generated $R$-submodules of the
fraction field $L$ of $R$.

\medskip
\noindent
{\bf Main \kern-1pt Result.\kern.5pt}~
{\it Every integrally closed fractionary\/ $R$-ideal is a multiplier~ideal.}

\medskip
After this paper was first submitted, we learned that independently of us C.\, Favre and M.\,
Jonsson had found a related proof \cite{FJ, \S6}. Their argument is given in the context of a
novel treatment of valuations of $R$. Though we thought initially that our proof applied only
to $\frak m$-primary ideals, Favre and Jonnson had no such restriction. This prompted us
to reexamine our proof, which we then~found could be made to apply to the general
case as well.

\goodbreak
\smallskip
Throughout, $(R,\frak m)$ and its fraction field~$L$ will be as above. ``Ideal" will
mean ``fractionary $R$-ideal." An ideal is {\it integral\/} if it is contained in $R$; and {\it
of finite colength\/} if it is either (integral and) $\frak m$-primary or the unit ideal.
 For brevity we use the classical term ``complete" instead of   ``\kern.5pt
integrally closed in
$L$."  Any 
complete ideal~$I$ is uniquely the product of a principal ideal
and a finite-colength complete ideal: $I=(I^{-1})^{-1}(II^{-1})$ where
the inverse of an ideal~$J$ 
is $J^{-1}\set\{\,x\in L\mid xJ\subset R\,\}$.

\subheading{1. Geometric formulation of the problem} The goal of this
section is to develop the geometric criterion Corollary 1.4.2
for an ideal to be a multiplier ideal, while laying the groundwork
for the proof in the next section that every complete ideal 
satisfies that criterion.

We begin by recalling
some preliminary definitions and known results. (For some historical
pointers to the development of the theory of complete ideals see the second
paragraph on the first page of \cite{L3}.)

For any complete
ideal $I$ there exists a {\it log resolution,\/} that is,  a proper birational
map $f:X\to \Spec(R)$ where $X$ is a regular scheme such that for any closed point 
$x\in X$ there exist $t_1$, $t_2$ generating the maximal ideal of $S\set\O_{X\<,\>x}$ 
such that $IS=t_1^{a_1}t_2^{a_2}S$ for some integers $a_1$,
$a_2$.  (In other words, $I\OX$ is invertible, with normal-crossing support.) To see quickly
that there is an~
$f\<$, composed of maps obtained by blowing up closed points, and such that the
$\OX$-ideal~$I\OX$ is at least {\it invertible\/}---assuming, as one clearly may, that $I$ is
$\frak m$-primary---one can use the Hoskin-Deligne formula \cite{L2, p.\,222, Thm.\, 3.1},
which shows that the length of the ``transform'' of~$I\/$ can be successively lowered by
suitable closed-point  blowups, until it vanishes, at which point $I\/$ generates an 
invertible\- ideal sheaf. Then well-known facts about embedded resolution of curves in
two-dimensional regular schemes ensure that with further closed-point blowups one can
reach the desired normal-crossing situation.

The group ${\Div(X)}$ of  {\it
$X\<\<$-divisors\/} is, by definition, the free
abelian group on the set of reduced irreducible one-dimensional subschemes of $X\<$. 
These subschemes, also called {\it prime divisors,} are of two kinds,  {\it affine\/}
resp.~{\it exceptional\/} according as
$f$ maps their generic point to a non-closed resp.~closed point of $\Spec(R)$. 
A~divisor~$D$ can be represented as a formal sum with integer coefficients,
\hbox{$D=\sum_Ed^{}_E E$}
where $E$ runs through all prime divisors and only finitely many of the integers~$d_E$ are
nonzero; an $E$ for which $d_E\ne0$ will be called  an {\it irreducible component\/} of
$D$. One has then the decomposition
$D=D\aff+D\exc$ where
$D\aff$ (resp.~$D\exc$), the {\it  affine\/} (resp.~{\it exceptional\/}) {\it part\/} of
$D$ is obtained by replacing all the~$d_E$ where $E$ is
exceptional (resp.~affine) by 0. We say $D$ is affine (resp.~exceptional) if 
$D=D\aff$ (resp.~$D=D\exc$). 
Exceptional divisors will also be called {\it $f\<\<$-divisors.}
They make up a subgroup $\Dive(X)\subset\Div(X)$.

Each prime divisor
$E$ gives rise to the discrete valuation
$v_{\<\<E}^{}$ whose valuation ring is the local ring
on~$X$ of the generic point of~$E$ (which may be assumed to be a subring
of $L$). The divisor $\blp t\brp$ of a nonzero $t\in L$ is
defined to be $\blp t\brp\set \sum_Ev_{\<\<E}^{}(t)E$. One has $t\in R\;\Leftrightarrow
\;v_{\<\<E}^{}(t)\ge 0$ for all affine
$E\;\Leftrightarrow\; \blp t\brp_\aff\ge0$.\vadjust{\kern1pt}

\penalty -2000
For $D=\sum_Ed^{}_E E$,  the sheaf $\OX(D)$ is the invertible $\OX$-module
sending any open $U\subset X$ to 
$$
\Gamma\bigl(U\<,\>\OX(D)\bigr):=\{\,t\in L \mid 
v^{}_E\<(t)\ge  -d^{}_E \text{ for all $E$ with generic point in } U\, \}.
$$
In particular,
$
\Gamma\bigl(X\<,\>\OX(-D)\bigr):=\{\,t\in L \mid 
\blp t\brp\ge D \}.
$
This is a complete ideal, integral iff $D\aff\ge 0$, and of finite colength  iff $D$ is
exceptional. Every invertible
$\OX$-submodule of the constant sheaf~$L$ is $\OX(D)$ for a unique $D$.
\vadjust{\kern 1pt}

The abelian group $\Div(X)$ is ordered, 
the positive (or {\it effective\/}) divisors being those $D$ such that
$d_E\ge0$ for all~$E$, or equivalently, $\OX(-D)\subset\OX$. An effective~$D\/$ can be
regarded as a one-dimensional subscheme of $X\<$,
with structure sheaf $\OD$ fitting in a natural exact sequence\vadjust{\kern-1\jot}
$$
0\to \OX(-D)\to \OX\to \OD\to 0.
$$\vskip-1\jot
\noindent If, moreover, $D$ is exceptional, then this subscheme
of~$X$ is projective over the field $R/\frak m$. In particular, the exceptional
prime divisors 
$E^1,E^2,\dots,E^s$  are isomorphic
to the projective line $\Bbb P^1_{\!R/\mkern-1mu\frak m}$; and any two of them
intersect\vadjust{\kern-1pt} transversally. (This can easily be shown by induction on the
number of blowups making up~$f\<$.)

For any invertible $\OX$-module $\L$ the {\it intersection
number\/} $\L\cdot E$ of $\L$ with an effective $f\<\<$-divisor $E$ is the degree of the
invertible $\OE$-module $\L_E\set \L\otimes\OE$:
$$
\L\cdot E =\chi_{\<E}^{} \O_{\mkern-2.5mu E}^{} - \chi_{\<E}^{}\> \L_E^{-1}
$$
where $\chi_{\<E}^{}\>\Cal M$ denotes the Euler characteristic of a coherent
$\OE$-module~$\Cal M$
($E$~being viewed as a projective curve over $R/\frak m$.)

For any $X\<\<$-divisor $F$ set 
\hbox{$F\<\<\cdot\<\< E\>\>\set \OX(F)\<\<\cdot\<\< E$}. This
intersection number extends uniquely to a  $\Bbb Z$-valued 
symmetric bilinear form on ${\Dive(X)}$ (see e.g., \cite{L1, \S13}).
If $F'=\blp t\brp+F$ for some $t\in L$ then $\OX(F')\cong\OX(F)$, and hence
$F'\<\<\cdot\<\< E=F\<\<\cdot\<\< E$.

An $X\<\<$-divisor $F$  is said to be numerically effective, 
{\it nef\/} for short,
if \hbox{$F\cdot E^i\ge 0$} for all $i=1,2,\dots,s\,$ 
$(\Rightarrow  F\cdot E\ge 0$ for all effective $E\in\Dive(X)$). $F$ is said
to~be {\it antinef\/} if $-F$ is nef. 

The following basic result is contained in \cite{L1, p.\,220,
Thm.\,(12.1)}. 

\proclaim{Theorem 1.1} 
An\/ $X\!$-\kern.5pt divisor\/ $D$ is nef if and only if\/ $\OX(D)$ is generated 
by its global sections.
\endproclaim

\proclaim{Corollary 1.1.1} 
If\/  $D$ is antinef and\/ $D\<\aff\ge 0$ then $D$ is effective.
\endproclaim

That's because, as above,  $\Gamma\bigl(X\<,\OX(-D)\bigr)\subset R$, and so by the
Theorem,
$$
\OX(-D)=\Gamma\bigl(X\<,\OX(-D)\bigr)\OX\subset\OX.
$$ 
(A simpler proof, not using Theorem 1.1, can be found in \cite{L1, p.\,238}.) 

\smallskip
If $I$ is a complete ideal and $I\OX$ is invertible then 
$I=\Gamma(X,I\OX)$, whence:
\proclaim{Corollary 1.1.2} {\rm(Cf.~\cite{L1, \S18}.)}
Sending\/ $E$   to\/ $\Gamma\bigl(X\<, \OX(-E)\bigr)$ defines an isomorphism 
from the $(\<\<$additive$)$ monoid of antinef\/ $\<X\<\<$-divisors
to~the 
$(\<$multiplicative\/$\<)$
monoid of those complete\/ 
ideals\/~$I$ such that\/ $I\OX$ is invertible.\footnotemark\  Under this isomorphism
antinef\/ $f\<\<$-divisors correspond to finite-colength complete ideals. 
\endproclaim
\footnotetext{It is a theorem of Zariski that a product of two complete
ideals is still complete \cite{ZS, p.\,385, Thm. 2$'$},
\cite{L1, p.\,209, Thm.\,(7.1)}.}

It is simple to show, by induction on the number of blowups making up~
$f$, that the intersection matrix $(E^i\cdot E^j)$ has determinant
$\pm 1$. Hence for each~$i$ there is a unique $f$-\kern.7pt divisor $G_i$ such
that $G_i\cdot E^j=0$ unless $j=i$, in which case $G_i\cdot
E^j=-1$; and for any $f$-\kern.7pt divisor $E$ it holds that 
$\>-E=\sum_{i=1}^s (E.E^i)G_i$.
Thus the monoid of antinef $f\<\<$-divisors is freely generated by
these~$G_i$. 
In other words, ``unique factorization'' holds in this
monoid---and therefore in the monoid of finite-colength complete ideals to which, by
Corollary 1.1.2, it is isomorphic.

An integral ideal $P\ne R$ is {\it simple\/} if whenever $P=IJ\>$ ($I$, $J$ integral ideals)
then either $I$ or $J$ is the unit ideal. For example, if $G_i$ is as above then 
the $\frak m$-primary ideal $P_i\set\Gamma\bigl(X\<,\OX(-G_i)\bigr)$ is simple, since if
$P_i$ is the product of $I$ and $J$ then it is also the product of the integral closures of
$I$ and $J\<$.

\proclaim{Corollary 1.1.3} {\rm (Zariski, \cite{ZS, p.\,386, Thm.\,3}.)} 
Every complete integral ideal is, in a unique way, the product of simple complete ideals.
\endproclaim

(As $R$
is a unique factorization domain, one reduces at once to the finite-colength
case. Also, it helps to
note that for any ideals $I$, $J$, 
if $IJ\OX$ is invertible then $I\OX$ and $J\OX$ are both invertible. )

\proclaim{Corollary 1.1.4} If\/ $I$\vadjust{\kern1pt} is a complete ideal with\/ 
$I\OX$ invertible, \hbox{$I\OX=\OX(-E),$} and\/
$P_i\set\Gamma\bigl(X\<,\OX(-G_i)\bigr)$ is the simple\/
$\frak m$-primary complete ideal corresponding to the above\/~$G_i\>,$ then\/
$II^{-1}=\prod_{i=1}^s P_i^{-E.E^i}\<\<.$ Thus $P_i$ divides\/
\hbox{$II^{-1}\Leftrightarrow E\cdot E^i\ne 0$.}
\endproclaim

Moreover, the valuations $v^{}_{\<\<E^i}$ associated to those $E^i$ such that 
$E\cdot E^i\ne 0$ are precisely the Rees valuations of $II^{-1}\>$ (i.e.,
those valuations whose valuation ring is the local ring of the
generic point of some
reduced  irreducible component of the closed fiber of the normalized blowup of
$I$). (See \cite{L4, p.\, 300, Prop.\,(4.4)}.)

\smallskip
The following Lemma%
\footnote{related to Enriques's ``principle of discharge''
 \cite{Z, p.\,28},} will be needed. 
\proclaim{Lemma 1.2}
Let\/ $E$ be an\/ $X\<\<$-divisor and\/ $I$ the complete\/ ideal\/
$\Gamma(X\<,\OX(-E))$. Then\/ $I\OX$ is invertible. Equivalently $($see $(1.1.2)),$ 
there is an antinef\/ \hbox{$E^-\ge E,$} 
such that\/ $I=\nmb\Gamma(X\<,\OX(-E^-));$ this\/ $E^-$ must be the\/
{\rm least} antinef divisor $\ge E.$
\endproclaim

\proof  Since $R$ is a unique factorization domain, there
exists a $t\in L$ such that $\blp t\brp_\aff=E_\aff\>$, and one can replace $(E,I)$ by
$(E-\blp t \brp,t^{-1}I$); so one may assume  that $E\/$ is exceptional, say
$E=\sum a_iE^i$, and that $I$ has finite colength. One may also assume $X\ne\Spec(R)$.
Then there exists an  antinef
$F\ge E$: pick~$n$ such that
$\frak m^n\subset I$, and define $F$ by $\frak m^n\OX=\OX(-F)$. Among all 
antinef $F=\sum b_iE^i\ge E$ choose one---call 
it $F_E$---for which
$\sum_i (b_i-a_i)$ has minimal value,  denoted $\sigma^{}_{\<\<E}$.

Procede by induction on $\sigma^{}_{\<\<E}$.
Suppose $\sigma^{}_{\<\<E}>0$ (otherwise there is nothing to prove), and that
the Lemma holds for all $E'$ with $\sigma^{}_{\<\<E'}<\sigma^{}_{\<\<E}$.
With $F\set F_E$  as above, there is an $i$ such
that $E\cdot E^i>0\ge F\cdot E^i$, and since 
$E^j\cdot E^i\ge0$ when $j\ne i$ (clearly) and
 $E^i\cdot E^i<0$,
\footnote{It is well-known, going back to Du Val, that the intersection matrix
$(E^i\cdot E^j)$ is negative definite, see e.g.,
\cite{L1, p.\,224, Lemma (14.1)}.}
 therefore $b_i>a_i\>$, whence
$F\ge E+E^i$. So $\sigma^{}_{\<\<E+E^i}<\sigma^{}_{\<\<E}\>$, and therefore
$\Gamma(X\<,\>\OX(-E-E^i))\OX$ is invertible. It suffices then to verify that
$
\Gamma\bigl(X\<,\>\OX(-E-E^i)\bigr)=I,
$
by applying the left-exact functor $\Gamma(X\<,-)$ to the 
natural exact sequence
$$
0\longrightarrow \OX(-E-E^i)\longrightarrow \OX(-E) 
\longrightarrow \OX(-E)\otimes \O_{\<\<E^i}
$$
and observing that since  $\OX(-E)\otimes\O_{\<\<E^i}$ 
has degree $-E\cdot E^i<0$, therefore
$$
\Gamma\bigl(X\<,\>\OX(-E)\otimes \O_{\<\<E^i}\<\bigr)=
\Gamma\bigl(E^i\<,\>\OX(-E)\otimes \O_{\<\<E^i}\<\bigr)=0.
\hbox to 0pt{$\hskip 40pt\square$\hss}
$$ 
\enddemo 

\nextpart{1.3} (Canonical divisors.) Let 
$Y\overset{g\,}\to{\to} X\overset{f\,}\to{\to}\Spec(R)$
be proper birational maps with $Y$ and $X$ regular schemes. 
By a theorem of Zariski and Abhyankar (see, e.g., \cite{L1, p.\,204,
Thm.\,(4.1)}) both $f$ and~$g$ are compositions of point blowups.

Let $F^1,F^2,\dots,F^t$  be the reduced irreducible components of 
$(fg)^{-1}\{\frak m\}$.
As before, the intersection matrix $(F^i\cdot F^j)$ has
determinant~$\pm 1$, and so there is a unique $fg$-divisor $K_g$ such that
$$
K_g\cdot F^i = 
\cases 
-F^i\<\<\cdot\<\< F^i\> -2 &\quad\text{if $g(F^i)$ is a point}\\
0 &\quad\text{otherwise.}
\endcases\tag{1.3.1}
$$
This $K_g$ is called the  {\it canonical divisor
of~$g$}.%
\footnote{One has $\O(K_{\<g})=H^0(g^!\OX)$ with $g^!$  as
in Grothendieck duality theory \cite{LS, p.\,206, (2.3)}.}
\vadjust{\kern1.5pt}%

For an $X\!$-\kern.5pt divisor~$D$, $g^*\<D$ denotes the $Y\<\<$-divisor 
whose coefficient at any prime divisor $F$ on $Y$ is $v^{}_{\<\<F}(\OX(-D)_{x})$, where
$x\in X$ is the
$g$-image  of the generic point of~$F\<$. There is a natural isomorphism
$g^*\OX(D)\iso\O_Y(g^*\<D)$.

\nextpart{1.3.2}  The following easily-checked properties characterize $K_{\<g}$ for
all~$g\>$:\vadjust{\kern1.5pt}\looseness=-1

$\bullet$ If $g$ is the blowup of a closed point $x\in X$ 
then $K_{\<g}=g^{-1}\{x\}$.

$\bullet$ If $Z\overset{h\,\>\>}\to{\to} Y\overset{g\,}\to{\to}
X\overset{f\,}\to{\to}\Spec(R)$ are proper birational maps with 
$Z$, $Y$ and~$X$ regular schemes, then
$$
K_{\<gh}= h^*\<K_{\<g} + K_{\<h}.
$$

\nextpart{1.4} (Multiplier ideals.)
For $D=\sum_E d_E E\in\Div(X)\otimes_{\Bbb Z}\Bbb R\>$ set 
$$
[D\>]\set\sum_E \,[d_E]E\in\Div(X)
$$
where $[d_E]$ is the greatest integer $\le d_E$.

\proclaim{Definition 1.4.1}
\rm{Let $I\/$ be a complete ideal,  $h:Y\to \Spec(R)$ a log
resolution of\/~$I$,
say $I\O_Y=\O_Y(-G)$, and let $c$ be a positive real number. 
The {\it multiplier ideal\/} $\J(I^c)$ is defined to be}
$$
\J(I^c)\set \Gamma\bigl(Y, \O(K_{\<h}-[c\>G\>])\bigr).
$$
{\rm Thus, by Lemma 1.2, }
$$
\J(I^c)\O_Y=\O_Y\bigl(-F\bigr)
$$
\rm{where $F\set\bigl([c\>G\>]-K_{\<h}\bigr)^-$ is the least 
antinef $h$-divisor $\ge[c\>G\>]-K_{\<h}$.}
\endproclaim
\penalty-1000
For the blowup $h_1:Y_1\to Y$ of a closed point $y\in Y$ one finds via (1.3.2) 
that the log resolution~$h$ can be replaced by the log resolution
$h\smcirc h_1$ without affecting~$\J(I^c)$. (When $c$ is not an integer, the
log-resolution hypothesis that
$G$ be a ``normal crossing divisor" is important.)
As any two log resolutions are dominated by a third, obtained from each of the two by a
sequence of point blowups,%
\footnote{By the above-mentioned theorem of Zariski and
Abhyankar, it suffices to principalize some ideal sheaf on one of the
log resolutions by a sequence of point blowups
(``elimination of indeterminacies''), which can be done e.g., via
the Hoskin-Deligne formula, as before.}
it follows that $\J(I^c)$ does not depend on the
choice of the log resolution $h$.

\proclaim{Corollary 1.4.2} A complete ideal\/ $J$ satisfies\/
$J=\Cal J(I^c)$ for some\/ $c,I$ iff for some log resolution\/ 
$h:Y\to\Spec(R)$ of\/ $J,$ say\/ $J\O_Y=\O_Y(-F),$ there is an 
antinef\/ $h$-divisor\/~$G$ and a real\/ $c>0$ such that 
$$
F=\bigl([c\>G\>]-K_{\<h}\bigr)^-\,.
\tag 1.4.2.1
$$
\endproclaim

\subheading{2. Proof of Main Result} Let $J$ be a complete ideal. To show that
$J=\Cal J(I^c)$ for some\/ $c,I$,
we will describe a log resolution $h:Y\to\Spec(R)$ of\/ $J\<$, 
and a $G$ and $c$ as in Corollary 1.4.2, such that if
$J\O_Y=\O_Y(-F)$ then (1.4.2.1) holds. (The number of suitable $(h, G,
c)$ will be enormous.)

Factor $JJ^{-1}$ as
$JJ^{-1}=\prod_{\ell=1}^{u} P_\ell^{e^{}_{\<\ell}}\ (P_\ell\text{ simple
complete, } e^{}_{\<\ell}>0)$---see Corollary 1.1.3 and the  paragraph
preceding it.
Let $f:X\to\Spec(R)$ be any log resolution of $J$, say
$J\OX=\OX(-F^{\>0})$. As in Corollary 1.1.4, one has  for each $\ell$,
$$
e^{}_{\<\ell}=-F^0\cdot E^\ell.\tag{2.0}
$$

We will first construct, for each
\hbox{$u$-tuple} \hbox{$N\set(n_1,n_2,\dots,n_u)$} of non-negative integers, a proper
birational map
\hbox{$g^{}_{\<N}:Y_N\to X$}, 
to be realized as a composition of closed point blowups, so that \hbox{$h_N\set f\smcirc
g^{}_{\<N}$} will be a log resolution of ~$J\<$.

For ease of expression we say ``blow up a closed point $x\in X$
generically, $n$~times'' to mean ``blow up $x_0\set x$ to get
$g_1:X_1\to X$, then blow up 
a closed point $x_1$ on $g_1^{-1}x_0$ but not on any other
irreducible component of $g_1^*F^0$ to get
$g_2:X_2\to X_1$, then blow up a closed point $x_2$ on $g_2^{-1}x_1$ but not
on \dots then blow up 
a closed point $x_{n-1}$ on $g_{n-1}^{-1}x_{n-2}$ but not on any other
irreducible component of  $(g_1\smcirc g_2\smcirc\cdots\smcirc
g_{n-1})^*\<F^0$ to get $g_n:X_n\to X_{n-1}$.''

Then with $g\set g_1\smcirc g_2\smcirc\cdots\smcirc g_n$ it
holds that:

\smallbreak
\noindent{\bf (2.1)} $g^{-1}x $ is a chain of $n$ integral
curves $D_{1},D_{2},\dots,D_{n}$ such that for $0<i<n$, 
$D_{i}\cdot D_{i+1}=1$ and  $D_{i}\cdot D_{i}=-2$, 
while $D_n\cdot D_n=-1$; and if $|j-i|>1$ then $D_{i}\cdot D_{j}=0$.
\vadjust{\kern1pt} 

(For the proof one can use, e.g., \cite{L1, p.\,229, middle, and p.\,
227, $\alpha)$ and $\beta)$}. Here, and subsequently, the reader may find it useful to do some rough sketches.)

As in Corollary 1.1.4, there corresponds to each $D_i$ an $\frak m$-primary simple
complete ideal
$Q_i\>$; and, we claim, {\it these $Q_i$ form a strictly decreasing
sequence\/} $Q_1> Q_2>\dots>Q_{\<n}\>$, with $Q_1$ strictly contained in each of
the simple ideals corresponding to the (one or two) prime $f$-divisors
$E^j$  passing through~$x$. 

Indeed,
let $G_{\<\<j}$ be the \hbox{$f$-\kern.7pt divisor} such that $G_{\<\<j}\cdot E^j=-1$
and $G_{\<\<j}\cdot E=0$ for every other prime $f$-divisor~$E$, and let 
\smash{$\widetilde Q_j\set\Gamma(X\<, \OX(-G_{\<\<j}))$}
be the corresponding simple complete ideal.  It follows from, e.g., \cite{L1,
p.\,227, $\alpha)$ and $\beta)$} that $g_1^*G_{\<\<j}$ is antinef; and the
corresponding simple complete ideal is
$$
\Gamma\bigl(X_1\>, \O_{\!X_1}(-g_1^*G_{\<\<j})\bigr)=
\Gamma\bigl(X\<, \OX(-G_{\<\<j})\bigr)=\widetilde Q_j\>.
$$
Further, with $E':=g_1^{-1}x$, 
let $G'$  be the $fg_1^{}$-\kern.7pt divisor 
such that $G'\cdot E'=-1$ and
$G'\cdot E''=0$ for every other prime $fg_1^{}$-\kern.7ptdivisor~$E''$. 
Since\vadjust{\kern.6pt} $g_1^*G_{\<\<j}+E'$ has intersection number 
$-1$ with $E'$ and $\ge 0$ with each~$E''$,\vadjust{\kern.6pt} therefore
$G'-g_1^*G_{\<\<j}-E'$ is antinef, hence effective (Corollary 1.1.1); and
consequently\vadjust{\kern.6pt}
$G'>g_1^*G_{\<\<j}$. Thus the simple complete ideal $\Gamma\bigl(X_1\>,
\O_{\!X_1}(-G')\bigr)$ is strictly contained in $\widetilde Q_j$.\vadjust{\kern.6pt} 

Continuing in this way we establish the claim.
 
\medbreak

 Now for each $\ell=1,2,\dots,u$, pick $e^{}_{\<\ell}$ distinct closed
points 
$x^{}_1,\dots,x_{e^{}_{\<\ell}}$ which
lie on $E^\ell$ but on no $E\ne E^\ell$ and blow up all of
these points generically,
$n^{}_{\<\ell}$ times. Then $Y_N$ is the resulting surface, and $g^{}_{\<N}$ is the
composition of all the blowups. It is easily seen that $(Y_N, g^{}_{\<N})$
does not depend (up to isomorphism) on the order in which the chosen
points are blown up---though that won't  be important.
\footnote{The initial $\sum_\ell e^{}_\ell$ 
points might be taken to be
the intersection of the
closed fiber on $X$ with a sufficiently generic curve 
$C$ in the linear system $|-F^0|$ (i.e.,
a divisor---having no component in the closed fiber---of 
the form $(j)-F^0$ with $j$ a sufficiently generic element of~$J$). Then at each stage
the point to be blown up could be taken to be a specialization of some nonclosed
point of~$C$.}

To simplify notation, fix $N$ and set $(Y,g)\set (Y_N, g^{}_{\<N})$ 
and $F\set g^*F^{\>0}$, so that $J\O_Y=\O_Y(-F)$. Also, 
set $h\set fg:Y\to \Spec(R)$.

For an $X\<\<$-\kern.7pt divisor $D$, we denote by $D^\#$ the {\it proper
transform\/} of $D$ on $Y\<$, obtained from  $g^*\<D=:\sum_E a_E E$ (where $E$ runs
through all prime divisors on~$Y$) by replacing  $a_E$ by 0 whenever $E$ is
$g\>$-{\it exceptional,} i.e., $g(E)$ is a closed point.  

For each $\ell=1,2,\dots,u$ and 
$x_{\!j^{}_{\<\ell}}\in E^\ell\ (j^{}_{\<\ell}=1,2,\dots, e^{}_\ell)$ let 
$$
\{\,E^\ell_{j^{}_{\<\ell}k^{}_{\<\ell}}
\mid k^{}_{\<\ell}=1,2,\dots,n^{}_{\<\ell} \,\}
$$
be the family of  prime $Y\<$-divisors 
whose $g$-image is $x_{\<j^{}_{\<\ell}}$, the
ordering of these curves by the index $k^{}_{\<\ell}$ conforming to
the ordering of the $D$'s in (2.1). These curves are all isomorphic to the projective line
$\Bbb P^1_{\!R/\mkern-1mu\frak m}$.

If $a^{}_\ell$ is the $E^\ell$-coefficient of the divisor
$F^{\>0}$, and $b^{}_\ell$ of the divisor $K_{\!f}$, then one finds
(using (1.3.2)) that
$$
\align
F=g^*\<F^{\>0}&= F^{\>0}{}^\# + 
 \sum_\ell \sum_{j^{}_{\<\ell}\<,\>k^{}_{\<\ell}}a^{}_\ell E^\ell_{j^{}_{\<\ell}k^{}_{\<\ell}},\\
g^*\<K_{\!f}&=K_{\!f}^\#+\sum_\ell
\sum_{j^{}_{\<\ell}\<,\>k^{}_{\<\ell}} b^{}_\ell E^\ell_{j^{}_{\<\ell}k^{}_{\<\ell}},\\
K_{\<g}&=\sum_\ell \sum_{j^{}_{\<\ell}\<,\>k^{}_{\<\ell}}  
  k^{}_\ell E^\ell_{j^{}_{\<\ell}k^{}_{\<\ell}}.
\endalign
$$

Set $G\set F+K_{\<g}$. Noting that
$F\cdot E^\ell_{j_{\<\ell}k^{}_{\<\ell}} =g^*\<F^0\cdot
E^\ell_{j_{\<\ell}k^{}_{\<\ell}}=0$ 
\hbox{\cite{L1, p.\,227, $\beta)$},}  and using (2.1) together with the preceding
expansion of  $K_{\<g}$ (or together with (1.3.1) and (1.3.2)),
one finds\vadjust{\kern.7pt} that for every $g$-exceptional prime divisor
$E$, $\>G\cdot E=0$ unless $E$ is one
of the curves 
$E^\ell_{j^{}_{\<\ell}n^{}_{\<\ell}}$ at the end of the chains
emanating from the $\sum_\ell e_{\<\ell}\>$ originally\vadjust{\kern.7pt}
chosen points (i.e., $\>g(E)$ is a point and $E\cdot E=-1$), in which
case $G\cdot E=-1$. Moreover,\vadjust{\kern.6pt} for any $f$-exceptional prime divisor
$D$,
$g^*\<F^0\cdot D^\#=F^0\cdot D$
\cite{L1, p.\,227, $\alpha)$}, and  $E^\ell_{j_{\<\ell}k^{}_{\<\ell}}\!\cdot D^\#=0$ if
$k_\ell>1$ (since then\vadjust{\kern-1pt} $E^\ell_{j_{\<\ell}k^{}_{\<\ell}}\!\cap D^\#$
is empty), so~using~(2.0) one finds\vadjust{\kern.7pt} that $G\cdot D^\#=0$. Thus
{\it $G$ is an antinef\/ $Y\!$-divisor.}

By Corollary 1.1.4, if $I$ is the corresponding complete ideal
$\Gamma\bigl(Y, \O_{Y}(-G)\bigr)$
then $II^{-1}$ is the product of the simple complete ideals corresponding to the
$\sum_\ell e^{}_{\<\ell}$  curves $E^\ell_{j^{}_{\<\ell}n^{}_{\<\ell}}$ having
self-intersection $-1$.

\smallskip

Here is a key technical point:
\proclaim{Lemma 2.2} For all sufficiently small\/ $\epsilon>0$ there
exists\/ $N$ such that 
$$
(1+\epsilon)G-K_{\<h}=F + A\qquad\bigl([A]\le 0\bigr)\tag 2.2.1
$$
where the coefficient of\/ $[A]$ at each\/ 
$E^\ell_{j^{}_{\<\ell}n^{}_{\<\ell}}$ \! and at each affine prime\/ $Y\!$-divisor is\/~
$\>0$. 
\endproclaim

\proof
Using (1.3.2), one transforms~(2.2.1) into the equality
$$
\epsilon(F+K_{\<g})-g^*\<K_{\!f}= A.
$$

More explicitly (see above)
$$
\epsilon\bigl(F^{\>0}{}^\#+ \sum_\ell  \sum_{j^{}_{\<\ell}\<,\>k^{}_{\<\ell}}  
\bigl(a^{}_\ell+k^{}_\ell\bigr)E^\ell_{j^{}_{\<\ell}k^{}_{\<\ell}}\bigr)
-\bigl(K_{\!f}^\#+\sum_\ell  
  \sum_{j^{}_{\<\ell}\<,\>k^{}_{\<\ell}}b^{}_{\<\ell}E^\ell_{j^{}_{\<\ell}k^{}_{\<\ell}}\bigr)=A.
$$
So to get (2.2.1) we can choose any $\epsilon>0$ such that the
coefficients of the \hbox{$X\<\<$-divisor} $\>\epsilon F^{\>0}-K_{\!f}$ are all $\><\<1$,
and then look for
$n^{}_{\<\ell}$ such that $\epsilon(a^{}_\ell+ k^{}_\ell)-b^{}_\ell<1$
for all $\ell\>$ and for all $k^{}_\ell\le n^{}_{\<\ell}\>$, while
$\epsilon(a^{}_\ell+n^{}_{\<\ell})-b^{}_\ell\ge 0$. These conditions mean
precisely that $n^{}_{\<\ell}$ satisfies the inequalities
$$
 1/\epsilon + b^{}_\ell/\epsilon -a^{}_\ell > n^{}_{\<\ell}
  \ge b^{}_\ell/\epsilon -a^{}_\ell
\quad (\ell=1,2,\dots,u).
$$ 
Clearly,  such integers $n^{}_{\<\ell}$ can be found if $\epsilon<1$.
\endproof

\goodbreak

For  $c=1+\epsilon$ and $N$ satisfying Lemma 2.2, and with
$h:Y\to \Spec(R)$ and $F$, $G$, 
as before, we have
$$
F'\set \bigl([c\>G\>]-K_{\<h}\bigr)^-\le F,
$$
so that 
$$
J'\set \Gamma\bigl(Y\<, \O(-F')\bigr)\supset \Gamma\bigl(Y\<,
\O(-F)\bigr)=J.
$$
Let us verify that $J'=J\  \bigl(=(J^{-1})^{-1}JJ^{-1}\bigr)$, thereby proving the main 
result.

Since\/ $G=F+K_g$ and\/ $F$ have the same affine part,  the
affine part of $[A]$ must  be\/ $0$, and hence  $F'$ and $F$ have the same affine part.
This means that $(J'{}^{-1})^{-1}=(J^{-1})^{-1}$. So $\>J'\<J'{}^{-1}\supset JJ^{-1}\<$,
and we need only show that these two $\frak m$-primary ideals are equal.

Recall that the valuations $v^{}_\ell\set v^{}_{\<\<E^\ell}$ are
just the Rees valuations of~$JJ^{-1}$. (See the remark following Corollary 1.1.4). So 
$$
JJ^{-1}=\{\,\rho\in R\mid v^{}_\ell(\rho) \ge v^{}_\ell(JJ^{-1}) \text{ for all
}\ell=1,2,\dots,u\,\}.
$$

Thus we need only show that for each $\ell$, the $E^\ell{}^\#\<$-coefficient
$a'_\ell$ of $F'$ is the same as that of $F$ (namely $a^{}_\ell$).
Let us say that $\ell$ is ``good'' if $a'_\ell=a^{}_\ell$ and ``bad''
if $a'_\ell<a^{}_\ell$.

If $\ell$  is good then since $F'\le F$ 
therefore 
$$
F'\cdot E^\ell{}^\#\le F\cdot E^\ell{}^\#=F^0\cdot E^\ell\,\overset{(2.0)}\to=-e_\ell\>. 
$$
Corollary~1.1.4 shows then that $J'\<J'{}^{-1}$ is divisible by $P_\ell^{e^{}_\ell}$.
\vadjust{\kern1pt}

Suppose $\ell$ is bad.
For $j\in[1,e^{}_{\<\ell}\>]$ and with $a'_{jk}$ 
the $E^\ell_{jk}$-coefficient of~$F'$\vadjust{\kern1pt}
it is easily seen that 
$
a'_\ell=:a'_{j 0}\le a'_{j 1}\le a'_{j 2}\le\dots\le
a'_{jn^{}_{\<\ell}}=a^{}_\ell.
$
\footnote{With $v\set v_{E^\ell_{j,\>k+\<1}}\!\!$ and $S\supset R$  the
regular local ring 
 blown\vadjust{\kern.8pt} up to give
$E^\ell_{j,\>k+1}$,\vadjust{\kern1pt} one finds:
\centerline{$a'_{j,\> k+1}=
v(J'S)=v\bigl(((J'S)^{-1})^{-1}\bigr)+v\bigl(J'S(J'S)^{-1}\bigr)=a'_{jk} +
v\bigl(J'S(J'S)^{-1}\bigr).$}}  
So\vadjust{\kern.8pt} there is a $k\in[1,n^{}_{\<\ell}\>]$
such that 
$
a'_{j,\>k-1}< a'_{jk}= a'_{j,\> k+1}=\dots=a^{}_\ell.
$
Then 
$$
F'\cdot E^\ell_{jk}=
\cases
a'_{j,\>k-1}-2a'_{jk}+a'_{j,\>k+1}<0\quad&\text{if }k<n^{}_{\<\ell}\>,\\
a'_{j,\>n^{}_{\<\ell}-1}-a'_{jn^{}_{\<\ell}}<0   \quad&\text{if }k=n^{}_{\<\ell}\>.
\endcases
$$
From Corollary 1.1.4 and the remarks after 2.1, one deduces that $J'\<J'{}^{-1}$
is divisible by a simple complete ideal $P'_{\ell j}<P_\ell$. This
being so for all $j$, and the $P'_{\ell j}$ being distinct (Corollary 1.1.2), 
it follows from Corollary 1.1.3 that $J'\<J'{}^{-1}$ is divisible by $P'_{\ell
1}P'_{\ell 2}\cdots P'_{\ell e^{}_\ell} < P_\ell^{e^{}_\ell}$.
Thus (by Corollary 1.1.3 again) the existence of a bad $\ell\>$ 
leads to a factorization of $J'\<J'{}^{-1}$ which contradicts $J'\<J'{}^{-1}\supset
JJ^{-1}$. 

So every $\ell$ is good, and $J'\<J'{}^{-1}=JJ^{-1}$.
$\quad\square$

\remark{Remarks} 1. By the choice of $\epsilon$, the
$E^\ell$-\kern.7pt coefficient of $\epsilon F^{\>0}-K_{\!f} $ is $\><\!1$, i.e.,
$\epsilon a^{}_\ell-b^{}_\ell<1$, i.e., 
$b^{}_\ell/\epsilon - a^{}_\ell>-1/\epsilon$. It could happen that  
$b^{}_\ell/\epsilon - a^{}_\ell< 0$ for all $\ell$. In this case one
can take $N=(0,0,\dots,0)$, and then $\>J=\J((1+\epsilon)J)$.\vadjust{\kern1pt}

2. The proof shows that if $J$ is a simple complete $\frak m$-primary ideal then there
is a simple complete $\frak m$-primary ideal $P\subset J$ and a $c>0$ such that
$J=\J(cP)$.\vadjust{\kern1pt}

3. By way of illustration of our method, let
$J$ be a simple
$\frak m$-primary ideal of order 16, whose successive transforms have orders
$(8,8,4,4,2,2,1,1,0,0\dots)$. (For the existence  of such a
$J$ see e.g., \cite{L4, p.\,298, Cor.\,(3.1)}.) Here $\ell=1$, and one calculates that
$a_1=426$, $b_1=46$. Then any $\epsilon\in(0,5/48)$ will do. Since
$46/(5/48)-426>15$, the least possible value of $N$ is 16, which is attained, e.g., when
$\epsilon=23/244$. One has then that $J=\J(P^{53/48})$, where $P\subset J$ is a simple
complete ideal of order~16 with successive transforms of orders 
$(8,8,4,4,2,2,1,1,...,1,0,0\dots)$ (18 ones). 

A simpler permissible choice of $\epsilon$ would be 1/12.  But then one would
need $N=126$ blowups before $P$ appeared.

There might be other methods of finding a representation  $J=\J(I^c)$ with
$I$ ``closer" to $J$ than here. But $I=J$ can never occur, because the inequality
$b_1<\epsilon a_1$ of example 1 could then be deduced. (Exercise.) On the other hand, 
this inequality does hold for any of the successive
transforms of the present $J$.\vadjust{\kern1pt}

4. Since the
$c=1+\epsilon$ we have considered can be arbitrarily close to 1, one may ask if it is
possible for $c$ actually to be 1. (This would be the case studied in \cite{L5}, where $\J(I)$
is called the {\it adjoint ideal\/} of $I$.) 

\goodbreak
For simple complete $\frak m$-primary $J\<$, the answer is given by:

\proclaim{Proposition 2.3} A simple complete $\frak m$-primary ideal\/ $J$ is of the
form\/~$\J(I)$ for some~\/ $I$ $\iff J\not\subset \frak m^2 \iff 
 J=(a,b^n)R$ for some integer\/ $n>0$ and\/ $a,b\in R$ such that\/ 
$(a,b)R=\frak m$.
\endproclaim

\proof
The last $\iff$ holds because $J\not\subset\frak m^2$ means
that $J$ contains an element~$a$ such that $R/aR$ is a discrete
valuation ring. Moreover, if $(a,b)R=\frak m$ and
$z\in R$ is integral over $J=(a,b^n)R$ then the canonical image of $z$ in
$R/aR$ is integral over---and hence is a multiple of---that of $b^n$,
whence $z\in J\<$, and thus $J$ is complete (and clearly simple). It
is an easy  exercise to show that for such a $J\<$, $J=\J(J^2)$. (One
could use \cite{L5, p.\,749, Prop.\,(3.1.2)}.)\vadjust{\kern1.5pt}

For a simple complete $\frak m$-primary $J\<$, let
$f:X\to\Spec(R)$ be a log resolution obtained by successively blowing up {\it base
points of~$J$}---closed points at which $J$~does not generate an invertible ideal---for as
long as such points are available. (As noted before, the Hoskin-Deligne formula guarantees
that this process terminates.) From \cite{L1, p.\,199, Prop.\,1.2 and 
p.\,203, Prop.\,3.1} it follows that this~$f$ is the minimal desingularization of the
blowup of $J$: for any log resolution $h$ of~$J$ there exists a map $g:Y\to X$
composed of point blowups such that
$h=fg$.

We have seen  before that there is a unique
exceptional prime $X\<\<$-divisor~$E\/$ such that 
$J\OX\cdot E\ne 0$. This $E\/$   satisfies
$E\cdot E=-1$:  for, there is at least one exceptional prime $X\<\<$-divisor~ $E'$ such
that
$E'\cdot E'=-1$, namely the closed fibre $f'{}^{-1}\{x'\}$
for the  blowup $f':X\to X'$ of $x'\in X'$ coming last
in the sequence of blowups composing to~$f$ (see \cite{L1, middle of p.\,229}); and 
$J'\set J\O_{\!X'\!,\>\>x'}$ is not invertible, from which one sees, with
$\frak m'$ the maximal ideal of $R'\set\O_{\!X'\<\!,\>\>x'},$ that
 $J'=d\>\frak m'{^s}$ for some $d\in R'$ and $s>0$, whence
$J\OX\cdot E'=s$, so that $E'=E$.

We claim that {\it if\/ $F\cdot F =-2\>$ for all exceptional prime divisors\/ $F\ne E$
 then\/} $J\not\subset\frak m^2\<$.\kern1.7pt%
(The converse is part of the exercise at the end of the first
paragraph.)
Indeed, this condition on the $F$'s means that among the base points of
$J$ no two are ``proximate" to the same one, and the conclusion follows 
from \cite{L4, p.\, 301, (3)}. \looseness=-1

Assume now that  $J=\J(I)$. From $J\/$ being  $\frak m$-primary it follows easily that so is
$I$. ($I$ and $\J(I)$ have the same gcd.) Let $h:Y\to\Spec(R)$ be a  log resolution of $I$,
obtained as above by blowing up base points of $I$, and say $I\O_Y=\O_Y(-G)$. Let $F$ be
an exceptional prime
$Y\<$-divisor. As above, we find that $F\cdot
F=-1\Rightarrow G\cdot F<0$. (We may assume that $F$~is the
closed fiber for the last blowup in some sequence of blowups composing to~$g$, because
if $g$ is the blowup of a point lying on an exceptional prime divisor~$F_1$ then 
$0>F_1\cdot F_1=g^*F_1\cdot F_1^\#=F_1^\#\cdot F_1^\# +1$.) Since
$K_{\<h}\cdot F=-F\cdot F -2$ and $G$~is antinef, we see that
 $G-K_{\<h}$ is antinef, and hence $\O_Y(K_{\<h}-G)=JO_Y$. 

So $J\O_Y$ is invertible, i.e., $h$ is a log resolution of $J$, and as above there exists a
$g:Y\to X$ composed of point blowups such that $h=fg$. Let
$F^1{}^\#, F^2{}^\#,\dots, F^n{}^\#$ be the proper transforms on~$Y$ of the
prime $X$-divisors $F^1, F^2,\dots, F^n$  
other than the above~$E$. 
Then $F^i{}^\#\<\<\cdot F^i{}^\#\le F^i\cdot F^i\le-2$\vadjust{\kern.6pt}
and so $K_{\<h}\cdot F^i{}^\#\ge 0$.
But by \cite{L1, p.\,227, $\beta)$},\vadjust{\kern-5pt}
$$
K_{\<h}\cdot F^i{}^\#\le (K_{\<h}-G)\cdot F^i{}^\#=JO_Y\cdot F^i{}^\#
=JO_X\cdot F^i=0,\vadjust{\kern-5pt}
$$
and thus $K_{\<h}\cdot F^i{}^\#=0$, i.e., $F^i{}^\#\<\<\cdot F^i{}^\#=-2$, whence, finally,
$F^i\cdot F^i =-2$. The above claim shows then that $J\not\subset\frak
m^2$.
\endproof 
\endremark

\vfill\eject

\enddocument